 \pdfoutput=1

\documentclass[twoside,11pt,a4paper,leqno]{article}
\usepackage{amssymb,amsmath,amscd,euscript,verbatim,array}
\usepackage[center,pagestyles]{titlesec}
\usepackage{geometry}
\usepackage{mathrsfs}
\usepackage{anysize}
\usepackage{multirow}
\usepackage{array}
\usepackage{tabularx}
\usepackage{fancyhdr}
\usepackage{indentfirst}
\usepackage{CJK}
\usepackage{graphicx}
\usepackage{color}
\usepackage{ifpdf}
\marginsize{3.5cm}{3.5cm}{2cm}{2cm}

\titleformat{\section}{\centering\normalsize}{\thesection.}{0.5em}{}
\titleformat{\subsection}{\normalsize\bfseries}{\thesubsection.}{0.5em}{}
\titleformat{\subsubsection}{\normalsize\bfseries}{\thesubsubsection.}{0.5em}{}
\newcommand{\N}{\mathbb{N}}
\newcommand{\Z}{\mathbb{Z}}
\newcommand{\R}{\mathbb{R}}

\newtheorem{Theorem}{Theorem}[section]
\newtheorem{Definition}[Theorem]{Definition}
\newtheorem{Lemma}[Theorem]{Lemma}
\newtheorem{Example}[Theorem]{Example}
\newtheorem{Proposition}[Theorem]{Proposition}
\newtheorem{Remark}[Theorem]{Remark}
\newtheorem{Corollary}[Theorem]{Corollary}

\newcommand{\gm}{\gamma}

\newcommand{\eps}{\varepsilon}

\newcommand{\T}{\mathbb{T}}
\newcommand{\bthm}{\begin{Theorem}}
\newcommand{\ethm}{\end{Theorem}}
\newcommand{\bpr}{\begin{Proposition}}
\newcommand{\epr}{\end{Proposition}}
\newcommand{\blm}{\begin{Lemma}}
\newcommand{\elm}{\end{Lemma}}
\newcommand{\bex}{\begin{Exercise}}
\newcommand{\eex}{\end{Exercise}}
\newcommand{\be}{\begin{equation}}
\newcommand{\ee}{\end{equation}}
\newcommand{\beal}{\begin{aligned}}
\newcommand{\enal}{\end{aligned}}
\newcommand{\brm}{\begin{Remark}}
\newcommand{\erm}{\end{Remark}}
\newcounter{item}[section]

\newcommand{\Proof}{\noindent\textbf{Proof}\hspace{0.3cm}}
\newcommand{\End}{\ensuremath{\hfill{\Box}}\\}
\renewcommand{\title}[1]{\begin{center}\textbf{\large #1}\end{center}}
\renewcommand{\author}[1]{\begin{center}\small #1\end{center}}
\renewcommand{\date}[1]{\begin{center}#1\end{center}}

\makeatletter \@addtoreset{equation}{section}
\makeatother
\newcommand{\PreserveBackslash}[1]{\let\temp=\\#1\let\\=\temp}
\newcolumntype{C}[1]{>{\PreserveBackslash\centering}p{#1}}
\newcolumntype{R}[1]{>{\PreserveBackslash\raggedleft}p{#1}}
\newcolumntype{L}[1]{>{\PreserveBackslash\raggedright}p{#1}}

\newcolumntype{I}{!{\vrule width 1pt}}
\newlength\savedwidth

\begin{document}
\vspace{10pt}
\title{\Large Quantitative Destruction and Persistence of Lagrangian Torus in Hamiltonian Systems}
\vspace{6pt}
\author{\sc Lin Wang}

\vspace{10pt} \thispagestyle{plain}
\begin{quote}
\small {\sc Abstract.}  For an integrable Hamiltonian systems with $d$ degrees of freedom ($d\geq 2$), we consider quantitatively the existence and non-existence of the flow-invariant Lagrangian torus with given frequency under the perturbation beyond the scope of the classical KAM method in the $C^r$ topology. As applications, the non-existence result  gives a partial answer to an open problem on non-existence of invariant circles by Mather \cite[p.212]{M4} from 1988. The existence result sheds a light on another open problem on the existence of invariant circles with lower regularity by Mather \cite[Problem 3.1.1]{MY} from 1998.
\end{quote}
\begin{quote}
\end{quote}
\begin{quote}
\end{quote} \vspace{25pt}

\section{\sc Introduction}

In 2013, C.-Q. Cheng and the author \cite{CW} constructed an example to show that:
{\it given an integrable Hamiltonian systems with $d$ degrees ($d\geq 2$) of freedom, each invariant Lagrangian torus can be destroyed by an arbitrarily small (bump) perturbation in the $C^{r}$ topology with $r<2d$. }

That result is almost sharp due to the fact proved by P\"{o}schel \cite{P} (see also Salamon \cite{Sa}):
{\it the KAM tori with constant type rotation vectors are preserved under an arbitrarily small perturbation in the $C^{r}$ topology with $r>2d$. }

From the physical point of view, it is more natural to consider a trigonometric
polynomial perturbation compared to the bump function used in \cite{CW}. Such kind of systems can be realized as a superposition of finite many harmonic oscillators. In this note, we revisit the converse KAM result in \cite{CW} from the quantitative perspective.  More precisely, we consider the following question.
\begin{itemize}
\item \textbf{Question 1:}  Let $P_N$ be a trigonometric
polynomial of degree $N$ and satisfy $\|P_N\|_{C^r}\leq \eps$. If $P_N$ destroys  the  Lagrangian torus with the rotation vector $\omega$ of an integrable Hamiltonian system, then {what are the relation among $\eps$, $N$, $r$ and the arithmetic property of $\omega$}?
\end{itemize}

An answer to Question 1 and its relevance to persistence result are collected in Theorem \ref{mathe} and Theorem \ref{excc} below. Compared to \cite{CW}, an obvious new difficulty in this note is that one can not localize the trigonometric polynomial perturbation, which behaves quite differently from the  bump function. Accordingly, more  dedicated  construction of the perturbation and careful analysis are involved. Besides, the perturbed Lagrangian torus shown by Theorem \ref{excc} can not be captured by the classical KAM method, since the perturbation there is too large beyond the scope of that method. Before giving precise statement of  the main results, we need to clarify the definitions of KAM torus and Lagrangian torus in Section 1.1, and introduce some arithmetics on the rotation vector in Section 1.2.

\subsection{KAM torus and Lagrangian torus}
We use, once and for all, $\T^d$ to denote a $d$-dimensional flat torus. Let $H: \text{T}^\ast\T^d\to\R$ be a $C^r$ ($r\geq 2$) Hamiltonian, and $\Phi_H^t: \text{T}^\ast\T^d\to \text{T}^\ast\T^d$ be its flow.
\begin{Definition}[KAM torus]\label{dd}$\bar{\mathcal
{T}}^d$ is called a $d$-dimensional KAM torus if
\begin{itemize}
\item $\bar{\mathcal {T}}^d$ is a Lipschitz graph over $\T^d$;
\item $\bar{\mathcal {T}}^d$ is invariant under the Hamiltonian flow
 $\Phi_H^t$;
\item there exists a diffeomorphism
$\phi:\ \T^d\rightarrow \bar{\mathcal {T}}^d$ such that
$\phi^{-1}\circ\Phi_H^t\circ\phi=R_\omega^t$ for any $t\in \R$,
where $R_\omega^t:\ x\mapsto x+\omega t$ and $\omega$ is called
the rotation vector (frequency) of $\bar{\mathcal {T}}^d$.
\end{itemize}
\end{Definition}

From the geometric point of view, a submanifold in a $2d$-dimensional symplectic manifold is called a
Lagrangian torus if it is homeomorphic to the torus $\T^d$ and the symplectic form vanishes on it.
In general, the rotation vector of certain orbit on an invariant torus is not well defined. See  \cite[XIII, Proposition 1.3]{H11} for a delicate example. Even though the rotation vector of each orbit exists,
they may not be the same. For instance, the Lagrangian torus may contain several orbits with different rotation
vectors. Various invariant Lagrangian graphs can be constructed easily by using Ma\~{n}\'{e}'s Lagrangian.

\begin{Example}\label{e1}
Denote the local coordinates of $\mathrm{T}^\ast\T^d$  by \[(x,y):=(x_1,\ldots,x_d,y_1,\ldots,y_d).\]  We consider the Hamiltonian $H:\mathrm{T}^*\T^d\to\R$ given by
\[H(x,y):=\frac{1}{2}\|y\|^2+\langle y,\chi(x)\rangle.\]
Here we fix, once and for all, $\|\cdot\|$ to denote the norm in the tangent and cotangent space induced by the flat metric on $\T^d$. This example is inspired by  Nemitsky and Stepanov \cite{NS}. The associated Lagrangian is so called Ma\~{n}\'{e}'s Lagrangian:
\[L(x,\dot{x}):=\frac{1}{2}\|\dot{x}-\chi(x)\|^2,\]
where the vector field is constructed as follows:
\[\chi(x):=\alpha \Psi(x),\]
where $\alpha=(\alpha_1,\ldots,\alpha_d)$ is a non-resonant vector, and $\Psi(x)$ is a $C^\infty$ function satisfying
\[\Psi(x_0)=0,\ 0<\Psi(x)\leq 1, \ \forall x\in \mathbb{T}^d\backslash \{x_0\}, \ \int_{\mathbb{T}^d}\frac{dx}{\Psi(x)}=+\infty.\]
By \cite[Section 6.16]{NS}, the invariant torus $\mathcal {T}^d$ coincides with the closure of the recurrent points, but $\delta_{x_0}$ is the unique ergodic measure.

Note that the constant function is a classical solution of the Hamilton-Jacobi equation $H(x,Du(x))=0$, where $Du:=(\frac{\partial u}{\partial x_1},\dots,\frac{\partial u}{\partial x_d})$. Then $\mathcal {T}^d$ coincides with 0-section of $\mathrm{T}^*\T^d$, which implies it is a Lagrangian graph. In fact, it is a $C^\infty$ graph based on the construction of $\chi(x)$.
By \cite[Theorem 2.13]{NW}, $\mathcal {T}^d$ is also the Aubry set in view of Aubry-Mather and weak KAM theory.
\end{Example}

From the dynamical point of view, the Lagrangian torus with the rotation vector $\omega$ is defined as follows.
\begin{Definition}[Lagrangian torus]\label{dd1}
$\mathcal {T}^d$ is a called a $d$-dimensional Lagrangian torus with the
rotation vector $\omega$ if
\begin{itemize}
\item $\mathcal {T}^d$ is a flow-invariant Lipschitz  Lagrangian graph over $\T^d$;
\item each orbit on $\mathcal {T}^d$ has the rotation vector $\omega$.
\end{itemize}
\end{Definition}

Here the rotation vector of the orbit on the Lagrangian torus is defined as (\ref{rotv}) below.
In view of \cite{ams}, the Lagrangian torus with resonant rotation vector (see Definition \ref{reson} below) is too fragile to be preserved unless the system is integrable. Correspondingly,  it is sufficient to consider how to destroy the Lagrangian torus with the
non-resonant rotation vector. To achieve that, we only need that the Lagrangian torus with the
non-resonant rotation vector satisfies the following condition:
\begin{itemize}
\item it is a flow-invariant Lipschitz  Lagrangian graph over $\T^d$;
\item there exists a dense orbit with  rotation vector $\omega$.
\end{itemize}
Equivalently, {\it it is a  flow-invariant, topologically transitive Lipschitz  Lagrangian graph over $\T^d$.}


By Herman \cite[Proposition 3.2]{H2}, a KAM torus with non-resonant rotation vector must be a Lagrangian torus.  Moreover, it can be represented by a graph of a differential. Note that the non-resonant rotation vector is essential to guarantee the tours being a Lagrangian submanifold. This phenomenon was first shown  by Arnaud (see \cite[Section 3.3]{H2}).
\begin{Example}\label{e133} Let $H:\T^2\times\R^2\to\R$ be a  Hamiltonian defined by
\[H(\theta_1,\theta_2,r_1,r_2)=\frac{1}{2}(r_1-\psi(\theta_2))^2+\frac{1}{2}r_2^2,\]
where $\psi:\T\to\R$ is a smooth non-constant function.
Then the torus $\{\theta_1,\theta_2,\psi(\theta_2),0\}$ is a non-Lagrangian invariant torus. Moreover, the restricted dynamics is
conjugated to a non-ergodic rotation of $\T^2$, i.e. the identity map.
\end{Example}

On the contrary, a flow-invariant Lipschitz Lagrangian torus may not be a KAM torus. If the rotation vector is a Liouvillean type (see Definition \ref{arr} below), those kinds of Lagrangian tori can be constructed by the celebrated AbC  (approximation by conjugation) method (see \cite{AK,FK,He6,KH}). If the rotation vector is a Diophantine type, one can refer to \cite{Arna2,AF} for the construction of non-differentiable invariant circles for area-preserving twist maps. By using Ma\~{n}\'{e}'s Lagrangian like Example \ref{e1}, one can obtain weakly mixing Lagrangian tori in light of \cite{Fa}. According to Definition \ref{dd}, all of the Lagrangian tori mentioned above are not KAM tori.

\subsection{Arithmetic properties of the rotation vector.}

To state our result precisely, we recall some arithmetic properties of the rotation vector.
\begin{Definition}[Resonant vector]\label{reson}
A vector $\omega\in\R^d$ is called resonant if there exists
$k\in\Z^d\backslash\{0\}$ such that $\langle \omega, k\rangle=0$. Otherwise, it is called
non-resonant.
\end{Definition}
For $
k:=(k_1,k_2,\ldots,k_d)$, we use $|\cdot|$ to denote
Euclidean norm,  i.e. $|k|=(\sum_{i=1}^d k_i^2)^{1/2}$.
By the Dirichlet approximation (\cite[Lemma 2.1]{C}), for any given vector $\omega\in \R^d$ with $d\geq
2$, there is a sequence of integer vectors $k_n\in\Z^d$ with
$|k_n|\rightarrow\infty$ as $n\rightarrow\infty$ such that
\begin{equation}\label{ap90}
\left|\langle\omega, k_n\rangle\right|<\frac{C}{|k_n|^{d-1}},
\end{equation}
where $C$ is a constant independent of $n$.

\begin{Definition}[Arithmetic properties]\label{arr}
\
\begin{itemize}
\item A
vector $\omega\in\R^d$ is called a $\tau$-well approximable if there exists a
positive constant $C$ as well infinitely many integer vectors
$k_n\in\Z^d$ such that
\begin{equation}\label{ap}
 |\langle \omega,k_n\rangle|\leq \frac{C}{|k_n|^{d-1+\tau}}.
\end{equation}
\item
A non-resonant $\omega$ is called a {Liouvillean type} if it is $\tau$-well approximable for all $\tau>0$. Otherwise, it is called a {Diophantine type}. Namely,  there exist constants $D>0$ and $\beta\geq d-1$ such that for all integer vectors
$k\in\Z^d\backslash\{0\}$,
\begin{equation}\label{apx}
 |\langle \omega,k\rangle|\geq \frac{D}{|k|^{\beta}}.
\end{equation}
Suppose $\omega$ is a Diophantine type. The infimum of the set of all $\beta$ for which (\ref{apx}) holds is called the Diophantine exponent of $\omega$.
\end{itemize}
\end{Definition}

\subsection{Statement of the main results}
 Given any $u\in C^r(\mathbb{T}^d)$. Let us recall the  $C^r$-norm in the sense of H\"{o}lder:
\[\|u\|_{C^r}:=\sup_{|\alpha|\leq [r]}\sup_{x\in \mathbb{T}^d}|D^{\alpha}u(x)|+\sum_{|\alpha|=[r]}\sup_{x,y\in\mathbb{T}^d}\frac{|D^\alpha u(x)-D^\alpha u(y)|}{|x-y|^{r-[r]}},\]
where $[r]$ is the integer part of $r$, $\alpha:=(\alpha_1,\ldots,\alpha_d)$, $|\alpha|:=\alpha_1+\cdots+\alpha_d$, and
\[D^\alpha u(x):=\frac{\partial^{|\alpha|}u}{\partial x_1^{\alpha_1}\cdots\partial x_d^{\alpha_d}}.\]  Given any $u\in C^\infty(\mathbb{T}^d)$, define
\[\|u\|_{C^\infty}:=\sum_{k\in \mathbb{N}}\frac{\arctan (\|u\|_{C^k})}{2^k},\]
which is not a norm. Endow the space $C^\infty(\mathbb{T}^d)$ with the translation invariant metric $\|u-v\|_{C^\infty}$. The $C^r$ (resp. $C^\infty$) topology is induced by the $C^r$-norm (resp. the translation invariant metric). For simplicity, we  use, once and for all,  $u\lesssim v$ (resp.
$u\gtrsim v$) to denote $u\leq Cv$ (resp. $u\geq Cv$) for some
positive constant $C$. We use $u\sim v$ to denote $\frac{1}{C}v\leq u\leq Cv$ for some
positive constant $C$.
The main result is stated as follows.
\begin{Theorem}\label{mathe}
Let $H_0:\mathrm{T}^\ast\T^d\to\R$ be an integrable Hamiltonian system  given by
\[H_0(y):=\frac{1}{2}\|y\|^2.\]
\begin{itemize}
\item  [(1)]
 Given $0<\eps\ll 1$,
let $\omega$ be a $\tau$-well approximable rotation vector. Then the Lagrangian torus with rotation vector $\omega$ of  $H_0$ can be destroyed by a trigonometric polynomial perturbation $P_N(x)$ satisfying $\|P_N\|_{C^{r}}\leq \epsilon$ with $r<2d+2\tau$.
\item [(2)] Given  $r\in [{0},2d+2\tau)$, let $\omega$ be a Diophantine type rotation vector with the  exponent $d-1+\tau$. Then there exists $P_N$ of degree {\small\[N=\mathcal{O}\left(\epsilon^{-\frac{d+\tau+1}{2d+2\tau-r}}\right)\]}with $\|P_N\|_{C^r}\sim \epsilon$ such that $H_0(y)+P_N(x)$ admits no Lagrangian torus with  rotation vector $\omega$.
\item [(3)] If $\omega$ is  a {Liouvillean type} rotation vector, then the Lagrangian torus with $\omega$ can be destroyed by  $P_N$ of degree $N=\mathcal{O}\left({{\epsilon}^{-1/2}}\right)$ satisfying $\|P_N\|_{C^\infty}\lesssim\epsilon$.
    \end{itemize}
\end{Theorem}

Some remarks on Theorem \ref{mathe} are listed as follows.
\begin{itemize}
\item {\bf On the optimality.} Item (1) of Theorem \ref{mathe} is optimal for $d=2$ by Herman's result \cite{H33} on the existence of the invariant circle with constant type frequency. It is still open if the KAM torus with certain frequency can be preserved under a small perturbation in the $C^{2d}$ topology for $d>2$. Recently, P\"{o}schel \cite{P2} obtained a KAM result with lower regularity for perturbed (not Hamiltonian) vector fields. We also refer the reader to \cite{A,TL} for recent progress on Hamiltonian systems from the KAM perspective.

 Regarding Item (2) and (3), we are not sure about the optimality on $N$. To verify that, we may need a more refined version of current KAM method. However, any result in this direction seems to be far from the present possibilities. It is also  related to the topic on the construction of particular examples with certain Lagrangian tori but not KAM ones. So far, there is still a series of open problems  (see \cite{FKO,MY} for instance).
\item {\bf Relation to previous works.} The author  obtained  some partial results in similar spirit as Theorem \ref{mathe} for  area-preserving twist maps \cite{W4} and symplectic twist maps \cite{W2}. Based on Item (1) of Theorem \ref{mathe}, the results in both   \cite{W2} and  \cite{W4} can be improved. A crucial ingredient  for the improvement is replacing the Cauchy estimate by the Bernstein inequality.  The result in \cite{W2} can also be improved in a similar way.  In view of the converse KAM result by Mather  \cite{Mm2} for the standard map, it seems that the invariant circles (or Lagrangian tori) may exist under a ``large" but ``special" perturbation. Theorem \ref{mathe} implies that in order to preserve a given Lagrangian torus, a trigonometric
polynomial of suitable large degree $N$ (compared to the $C^\infty$ perturbation) is not special enough to relax the size of the perturbation in the $C^r$ topology.

 In the 1980's, Herman \cite{H1} (resp. Mather \cite{M4}) proved a remarkable result on destruction of the invariant circle with Diophantine frequency in the $C^r$ ($r<3$) topology (resp.  Liouvillean frequency  in the $C^\infty$ topology) for area-preserving twist maps.   Theorem \ref{mathe}  provides  quantitatively higher dimensional generalizations of their results.
We refer the reader to \cite{Fo} (resp. \cite{B2}) for the issue on destruction of invariant circles (resp. Lagrangian tori) in the real-analytic topology.
\end{itemize}

Note that the construction we give in  Theorem \ref{mathe} is a classical mechanical system. The Maupertuis principle allows us to regard classical trajectories as reparametrized geodesics of the Jacobi-Maupertuis metric on configuration space.

\begin{Corollary}\label{mathxc}
Let $\T^d$ be a flat Riemannian $d$-torus.
 Given  a  $\tau$-well approximable rotation vector $\omega$, one can find an arbitrarily small trigonometric polynomial perturbation  in the $C^r$ neighborhood with $r<2d+2\tau$ of the flat metric such that
there is no dense minimal geodesics with the rotation vector $\omega$.
\end{Corollary}

See \cite{B} for the definition of the minimal geodesics. This kind of geodesics was first studied by Morse \cite{Mor} in which they are called ``Class A".

In general, there are several examples of invariant Lagrangian tori which can not be captured by the KAM approach (see \cite{AK,Arna2,AF,FK,He6,KH} for instance). However, little  is known about the existence of the Lagrangian torus beyond the KAM approach if we require that the torus has a Diophantine type rotation vector and the system is a $C^r$ perturbation of the integrable one.
More precisely, one ask the following
\begin{itemize}
\item \textbf{Question 2:} Given $0<\eps\ll 1$,  $0\leq r<2d+2\tau$, let $\rho$ be a rotation vector with the Diophantine exponent $d-1+\tau$ and let $H$ be a Hamiltonian defined by
\[H(x,y)=\frac{1}{2}\|y\|^2+V(x),\quad  (x,y)\in \T^d\times\R^d.\]
Assume $\|V\|_{C^{r}}<\eps$. If the flow generated by $H$ admits a $C^1$ Lagrangian graph with rotation vector $\rho$, can we find $s_0>2d+2\tau$ such that $\|V\|_{C^{s_0}}<\eps$?
\end{itemize}
Here  we denote the rotation vector by $\rho$ instead of $\omega$ in order to distinguish from the notation of $C^\omega$ topology. 

We will give a negative answer to Question 2. 
As preparations, let us recall the definition of $C^\omega$ topology in the following.
\begin{Definition}[$C^{\omega}$ topology]
Let $\mathcal{A}$ be the set of all real analytic function $\varphi:\R^d\to\R$. For $r>0$, let $\mathcal{A}_r$ be the subset of $\varphi\in \mathcal{A}$ which admits a bounded holomorphic extension $\Phi$ in the strip
\[S_r:=\{z:=(z_1,\dots,z_d)\in \mathbb{C}^d\ |\  \|\mathrm{Im}z\|:=\max\{|\mathrm{Im}z_1|,\ldots,|\mathrm{Im}z_d|\}\leq r\}.\]
Denote the norm
\[\|\varphi\|_r:=\max_{S_r}|\Phi|.\]
Then $\mathcal{A}_r$, equipped with the norm $\|\cdot\|_r$ is a Banach space. For $r>r'$, $\mathcal{A}_r$ embeds as a subspace of $\mathcal{A}_{r'}$, and $\mathcal{A}_r$ embeds as a subspace of $\mathcal{A}$ for all $r>0$.
The $C^\omega$ topology on $\mathcal{A}$ is introduced as the direct limit of $\mathcal{A}_r$ as $r\to 0^+$. It can be described by the following fundamental system of neighborhoods of zero function. For any strictly positive real function $\eps:\R^+\to\R^+$, we let
\[\mathcal{U}_\eps:=\{\varphi\in \mathcal{A}\ |\ \exists r>0 \mathrm{\ such\ that\ }\varphi\in \mathcal{A}_r\ \mathrm{with}\ \|\varphi\|_r<\eps(r)\}.\]
\end{Definition}

  As an application of Theorem \ref{mathe}, we have
\begin{Theorem}\label{excc}
Let
\[H_0(y):=\frac{1}{2}\|y\|^2.\]
 Let $\rho$ be a Diophantine type rotation vector with the  exponent $d-1+\tau$.  Given $\eps>0$ and  $r\in [0,2d+2\tau)$,  there exists a  potential $V:\T^d\to\R$ of class $C^\omega$  such that
 \begin{itemize}
 \item  [(1)]
$\|V\|_{C^r}<\eps$;
\item [(2)]
$\|V\|_{C^s}\geq \eps$ for all $s> 2d+2\tau$;
 \item [(3)]
 $H_0(y)+V(x)$ has the unique Lagrangian graph $\mathcal {T}^d$ with  rotation vector $\rho$,
and the dynamics on $\mathcal {T}^d$ is at least Lipschitz conjugated to the linear flow on $\mathbb{T}^d$;
\item [(4)] $\mathcal {T}^d$ is a  graph of class $C^1$.
\end{itemize}
\end{Theorem}
Thanks to \cite[Theorem 2]{Arna0}, the $C^1$ smoothness of the graph directly follows from the dynamics on $\mathbb{T}^d$. Namely, Item (4) follows from Item (3).
In view of Item (2), the perturbed Lagrangian torus in Theorem \ref{excc} can not be captured by the classical KAM method. By Herman \cite{H33}, Item (2) can be improved to $s\geq 4$ for $d=2$ and $\tau=0$ (i.e. constant type rotation number).
\subsection{On two problems by Mather}
\subsubsection{A problem from 1988}
As an application of Theorem \ref{mathe}, we give a partial answer to an open problem by Mather \cite[p.212]{M4} from 1988. To keep consistency with the notations used in \cite{M4}, we denote $f$ to be a diffeomorphism of $\R^2$ denoted by $f(x,y)=(X(x,y),Y(x,y))$. Let $f$ satisfy:
\begin{itemize}
\item {\it periodicity:} $f\circ T=T\circ f$ for the translation $T(x,y)=(x+1,y)$;
\item {\it twist condition:} the map $\psi:(x,y)\mapsto(x,X(x,y))$ is a diffeomorphism of $\R^2$;
\item {\it exact symplectic:} there exists a real valued function $h$ on $\R^2$ with $h(x+1,y)=h(x,y)$ such that
    \[YdX-ydx=dh.\]
\end{itemize}
Then $f$ induces a map on the cylinder denoted by $\bar{f}$: $\T\times\R\mapsto \T\times\R$ ($\T=\R/\Z$). $\bar{f}$ is called an exact
area-preserving  twist map.  Let $\mathscr{T}^{\infty}$ be the set of $\bar{f}\in C^\infty$. Let $\tilde{\mathscr{T}}^{\infty}$ be the set of  $f\in C^\infty$.

\vspace{1ex}

In \cite[p.212]{M4}, Mather raised the following problem: {\it Suppose $\omega$ is a Diophantine number. The infimum of the set of all $N$ for which there exists $C>0$ such
that $|q\omega-p|>C|q|^{-N}$ for all $q, p\in \Z\backslash \{0\}$ is called the Diophantine exponent of $\omega$.   \cite[Proof of Theorem 2.1]{M4} shows that for every positive integer $r$ there is a
number $\gm$ such that if  $\bar{f}\in \mathscr{T}^{\infty}$, and $\omega$ has Diophantine
exponent $\geq \gm$, then in any $C^r$
 neighbourhood of $f$ in $\tilde{\mathscr{T}}^{\infty}$,
there is a $g$ such that $g$
has no homotopically non-trivial invariant circle $\Gamma$ such that $\rho(g,\Gamma)=\omega$. (Here $\bar{g}$
is the unique element of $\mathscr{T}^{\infty}$ of which $g$ is the lift to the universal cover.) Let $\gm(r)$ be the infimum of all such $\gm$. Let $\alpha$ be the infimum of all positive numbers $\alpha$ such
that $\gm(r)=\mathcal{O}(r^\alpha)$ as $r\to\infty$. \cite[Proof of Theorem 2.1]{M4} shows that $\alpha\leq 2$. On
the other hand, KAM theory shows that $\alpha\geq 1$ (see e.g. Salamon \cite{Sa}). {Find $\alpha$.}}

\vspace{1ex}

In light of Theorem \ref{mathe}, we formulate an alternative version of Mather's problem for Hamiltonian systems: {\it Let
\[H_0(y):=\frac{1}{2}\|y\|^2.\]
Assume $\omega$ is a Diophantine type, and $\omega$ has the
exponent $\geq \gm$. Suppose the Lagrangian torus with the rotation vector $\omega$ is destroyed by an arbitrarily small trigonometric polynomial perturbation in the $C^r$ neighbourhood of $H_0$. Let $\gm(r)$ be the infimum of all such $\gm$. Let $\alpha$ be the infimum of all positive numbers $\alpha$ such
that $\gm(r)=\mathcal{O}(r^\alpha)$ as $r\to\infty$. { Find $\alpha$.}}

 By Theorem \ref{mathe}(1), for a given  integer
\[r\in \left[2d+2\tau-\frac{1}{2},2d+2\tau+\frac{1}{2}\right),\]
where $d\geq 2$, $\tau\geq 0$,
we have
\[\gm(r)\leq d+\tau.\]
On the other hand, following \cite{P} (see also \cite{Sa}), we have
\begin{Proposition}\label{KK1}
Let $\omega$ be a rotation vector with the Diophantine exponent $d+\tau-2$. The KAM torus with $\omega$  is preserved under an arbitrarily small perturbation  in the $C^{2d+2\tau-2+\eps}$  neighbourhood of $H_0$.
\end{Proposition}

Note that for $0<\eps\ll 1$,
\[2d+2\tau-2+\eps<2d+2\tau-\frac{1}{2}\leq r.\]
 It follows that
\[\gm(r)\geq d+\tau-2.\]
It is clear to see  that $r\to\infty$ is equivalent to $\tau\to\infty$. Then
\[\lim_{r\to\infty}\frac{\gm(r)}{r}=\frac{1}{2},\]
which implies $\alpha=1$. Note that $\alpha$ is independent of the degree of freedom of the Hamiltonian system. 

\begin{Remark}
If we consider a nearly integrable area-preserving twist map of class $C^\infty$, the result $\alpha=1$  can be obtained following Herman's work \cite[Chapter II]{H1}. For the $C^\infty$ Hamiltonian system with multi-degree of freedom, one can also get  $\alpha=1$ following the work by Cheng and the author \cite{CW}, although this issue was not stated in both aforementioned works. The content here in this note shows that the answer to Mather's problem  is the same even if we use more special trigonometric polynomial perturbations.  Meanwhile, this problem is still open if we consider a general $C^\infty$ area-preserving twist map without any nearly integrable assumption.
\end{Remark}

In the case of twist maps, it is known by Birkhoff that each invariant circle must be Lipschitz. Then each (homotopically non-trivial) invariant circle is a Lagrangian submanifold in the sense that the tangent space is defined almost everywhere. Based on \cite{H11}, the dynamics on the invariant circle with irrational rotation number is much simpler than the higher dimensional case. To be more precise, let $\Gamma$ be the invariant circle with irrational rotation number $\omega$ that is preserved by $\bar{f}$. Then the following statements are equivalent:
\begin{itemize}
\item [(1)] $\bar{f}|_{\Gamma}$ is $C^0$ conjugated to $R_\omega$;
\item [(2)] there exists a strictly ergodic  probability measure $\mu$ with rotation number $\rho(\mu)=\omega$;
\item [(3)] each orbit is dense, and has rotation number $\omega$;
\item [(4)] there exists a dense orbit with rotation number $\omega$.
\end{itemize}
Consequently, the dynamics on a flow-invariant, topologically minimal,  Lipschitz Lagrangian graph over $\T^1$ must be $C^0$ conjugated to $R_\omega$. However, this kind of invariant circle may not be a KAM circle on which  the dynamics is at least $C^1$ conjugated to $R_\omega$. See a counterexample constructed by Avila and Fayad \cite{AF}.

\subsubsection{A problem from 1998}
Based on a similar idea in the proof of Theorem \ref{excc}, we give a new evidence to an open problem by Mather \cite[Problem 3.1.1]{MY} in 1998.

{\it Does there exist an example of a $C^r$  area-preserving twist map with an invariant curve which is not $C^1$ and that contains no periodic point? (separate the
question for each $r\in [1,+\infty]\cup \{\omega\}$)

}

We refer the reader to \cite{Arna2,AF} for recent progress on the existence of $C^1$ examples. The problem for $r\geq 2$ was raised again by Fayad and Krikorian in ICM2018 (see \cite[Question 26]{FKO}). In light of Mather's problem from 1998, one may ask an alternative question with a similar spirit.

\begin{itemize}
\item \textbf{Question 3:} Given a $C^\infty$ (reps. $C^\omega$)  area-preserving twist map with an irrational frequency invariant curve which is formed as a graph of $\psi$, what is the { minimal regularity} of   $\psi$?
\end{itemize}

Based on Herman-Yoccoz's global theory for circle diffeomorphisms, we have
\begin{Theorem}\label{neev11}
Given $0<\eps\ll 1$, let $\rho$ be a constant type  frequency. Let $f_0:(x,y)\mapsto (x+y,y)$ be an integrable twist map. Then there exists a sequence of $C^\infty$
area-preserving  twist maps $\{f_n\}_{n\in \mathbb{N}}$ satisfying
 \[\|f_n-f_0\|_{C^{3-\eps}}\to 0,\quad\|f_n-f_0\|_{C^3}\nrightarrow 0, \]
such that $f_n$ has an invariant circle $\Gamma_n:=\{(x,\psi_n(x))\ |\ x\in \mathbb{T}\}$ with frequency $\rho$. Moreover, $\psi_n$ is not of class $C^\infty$.
\end{Theorem}

\begin{Remark}
The proof of Theorem \ref{neev11} is essentially inspired by \cite[Corollary 3.5]{H1}. Honestly, we are not sure if we replace $C^\infty$ by $C^\omega$. The main difficulty comes from the fact that $C^\omega(\T)$ is not  a complete metric space such that we can not obtain the openness of certain set following from the inverse mapping theorem (see \cite[Theorem 2.6.1]{H1}).
\end{Remark}

Based on \cite{Sa}, we know that for a $C^\infty$ area-preserving twist map, if the invariant circle with constant type frequency is graph of class $C^r$ with $r>4$, then it is of class $C^\infty$. In \cite{KO},  Katznelson and Ornstein proved that if $f$ is a $C^{3+\gm}$ ($\gm>0$) area-preserving surface diffeomorphism that admits a $C^{2+\eps}$ invariant curve  $\psi$ with a ``good"  frequency, then $\psi$  is of class $C^{2+\gm'}$ for all $\gm'<\gm$. Here a ``good" frequency $\alpha$ means its continued fraction expansion $[a_1,a_2,\ldots]$ satisfies $a_n=\mathcal{O}(n^2)$, which is more general than the constant type frequency. Combining \cite{Sa} and \cite{KO}, we have

\begin{Theorem}\label{neev}
There exists a $C^\infty$ area-preserving twist map with  constant type frequency invariant  curve $\Gamma:=\{(x,\psi(x))\ |\ x\in \mathbb{T}\}$, and  $\psi$ is at most of class $C^2$.
\end{Theorem}

\begin{Remark}
With a numerical method, it was shown by Olvera and Petrov \cite{OP} that for some  area-preserving twist maps (for example, the standard map), the critical invariant circle (i.e. its rotation number equal to the golden ratio $\frac{\sqrt{5}-1}{2}$) is a $C^{1+\kappa}$ graph where $\kappa\in (0.7,1)$. The dynamics on the graph is conjugated to the rigid rotation by a $C^{\kappa'}$ function with  $\kappa'\in (0.7,1)$. This numerical evidence gives rise to further difficulty to reduce the regularity of $\psi$ in Theorem \ref{neev}.
\end{Remark}

The lower regularity of $\psi$ obtained in Theorem \ref{neev} is essentially based on Herman-Yoccoz's global theory for circle diffeomorphisms. Due to the lack of the corresponding theory for higher dimensional cases, one can not conclude the lower regularity of  the Lagrangian graph in Theorem \ref{excc}.
\subsection{Strategy of the proof of Theorem \ref{mathe}}
Inspired by Bessi \cite{B2}, we  use the variational method developed by Mather \cite{M5}. An autonomous Hamiltonian is called a Tonelli Hamiltonian if it satisfies strictly convexity and superlinearity with respect to momentum. Correspondingly, a Tonelli Lagrangian can be  defined by the Legendre transformation.   For Tonelli
Hamiltonian systems, a flow-invariant Lagrangian torus is  the graph over $\T^d$ if  it is Hamiltonianly isotopic to the zero-section (see \cite{Arna,BP}). By Herman \cite[Theorem 8.14]{H2}, a flow-invariant $C^0$ Lagrangian graph must be Lipschitz continuous. Here $\mathcal {T}^d$ is called  a $C^0$ Lagrangian graph if
 \[\mathcal {T}^d=(x,c+D\eta(x)),\]
 where  $\eta:\T^d\to\R$ is a $C^1$ function,
$c\in \R^d$ is a constant vector and $D\eta:=(\frac{\partial \eta}{\partial x_1},\dots,\frac{\partial \eta}{\partial x_d})$.  It was  shown by \cite{FGS} that the KAM torus with a given non-resonant rotation vector is unique in this setting.
\begin{Definition}[Action minimizing curve]
An absolutely continuous curve $\gamma\in C^{ac}([t_1,t_2],\T^d)$
is called an action minimizing curve  if
\[\int_{t_1}^{t_2}L(\gamma(t),\dot{\gamma}(t))dt\leq \int_{t_1}^{t_2}L(\bar{\gamma}(t),\dot{\bar{\gamma}}(t))dt,\]
for  every $\bar{\gamma}\in C^{ac}([t_1,t_2],\T^d)$ satisfying
\begin{itemize}
\item $\gamma(t_1)=\bar{\gamma}(t_1)$, $\gamma(t_2)=\bar{\gamma}(t_2)$;
\item $\gamma$ and $\bar{\gamma}$  are in the same homotopy
class. Equivalently, their lifts to $\R^d$ connect the same points.
\end{itemize}
\end{Definition}
Based on \cite{M5}, the action minimizing curves always exist and satisfy the Euler-Lagrange equation generated by $L$. For Tonelli Hamiltonians, a classical result by Herman \cite{H2} asserts that each orbit on the invariant Lipschitz Lagrangian graph is an action minimizing
curve.
To destroy this kind of  torus with non-resonant rotation vector, it is sufficient to prove that there exists a  neighborhood of certain point on the Lagrangian torus such that no action minimizing curve passes through.

\subsection{Strategy of the proof of Theorem \ref{excc}}

Let $\rho$ be a Diophantine type rotation vector with the  exponent $d-1+\tau$. Note that the perturbation $V_n$ is not small in the $C^s$ topology with $s>2d+2\tau$. Consequently, the proof of the existence of Lagrangian torus  is beyond the scope of the classical KAM method under such kind of conditions. Inspired by Herman \cite{H1}, we use a topological argument to prove this theorem. Let $\mathcal{V}$ be  the set of $C^\omega$ potentials $V:\mathbb{T}^d\to\R$. We choose $\mathcal{U}^r_\delta\subseteq \mathcal{V}$ to be an open and connected $\delta$-neighborhood of $V\equiv 0$ in the $C^r$ topology with $r<2d+2\tau$. Let $\mathcal{U}({\rho},r,\delta)\subseteq \mathcal{U}^r_\delta$ to be an open neighborhood  of $V\equiv 0$ in the $C^\omega$ topology.  For $\mathcal{U}({\rho},r,\delta)$ and $\rho$ introduced above, by using the KAM method developed by Salamon and Zehnder \cite{SaZ} (see also \cite{Sa}), we conclude that there exists $\delta>0$ such that for all $s>2d+2\tau$, if  $V\in \mathcal{U}({\rho},s,\delta)$, then the flow generated by the Hamiltonian $H_0+V$ admits the KAM torus with frequency $\rho$. Accordingly,  we consider the set $W_\rho:=\cup_{s>2d+2\tau}\mathcal{U}({\rho},s,\delta)$ and its closure  $\overline{W_\rho}$ in the $C^\omega$ topology. Based on a compactness argument, we know that for each potential $V$ in $\overline{W_\rho}$, the Hamiltonian flow generated by $H_0+V$ still has a Lagrangian torus with  frequency $\rho$. Moreover, Theorem \ref{mathe} implies
\[\Delta({\rho},r,\delta):=\overline{W_\rho}\backslash W_\rho\neq \emptyset.\]
Then each element in $\Delta({\rho},r,\delta)$ is qualified to be a potential for verifying Theorem \ref{excc}. The uniqueness of the Lagrangian torus with frequency $\rho$ is guaranteed by \cite{FGS}.

The note is organized as follows. In Section \ref{2s}, we introduce a change of coordinates to reduce the Lagrangian torus with a general rotation vector $\omega$ to certain torus with ``controlled" first two components of $\omega$. In Section \ref{3s}, we give a construction of the trigonometric polynomial perturbation by using an enhanced version of  Jackson's approximation. In Section \ref{41s}, we reduce Theorem \ref{mathe} to Proposition \ref{keyy}, and recall some facts about the estimates on the action of the pendulum and the speed of the action minimizing orbit. In Section \ref{4s}, we complete the proof of Proposition \ref{keyy}. The proofs of Theorem \ref{excc} and Theorem \ref{neev11} are provided in Section \ref{66} and Section \ref{77} respectively. This work is inspired by  \cite{CW}.  Section \ref{2s} and \ref{41s} in this note are parallel to  \cite{CW} essentially.

\section{\sc The change of coordinates}\label{2s}

The system we consider consists of one pendulum, a rotator with
$d-1$ degrees of freedom and a perturbation coupling of them.

We are concerned with a sequence of Hamiltonian functions
\begin{equation}\label{ha}
H_n(x,y)=\frac{1}{2}\|y\|^2-P_{N}(x),
\end{equation}
where $N$ (depending on $n$) denotes the degree of the trigonometric polynomial $P_N(x)$.
Since $H_n$ is quadratic with respect to $y$, by the Legendre
transformation, the Lagrangian function associated to $H_n$ is
\begin{equation}\label{Lo}
L_n(x,\dot{x})=\frac{1}{2}\|\dot{x}\|^2+P_{N}(x),
\end{equation}
where $\dot{x}=\frac{\partial H_0}{\partial y}$.    The following Lemma was proved in \cite{CW}.
\begin{Lemma}\label{kkk}
For any $k_n\in\Z^d\backslash\{0\}$ satisfying (\ref{ap}), there exists an integer vector $k'_n\in \Z^d$ such that
\begin{equation}\label{kn1}
\langle k_n,k'_n\rangle=0,\quad|k'_n|\sim |k_n|\quad\text{and}\quad
|\langle k'_n,\omega\rangle|\sim |k_n|.
\end{equation}
\end{Lemma}

 We choose two sequences of $\{k_n\}_{n\in \N}$ and $\{k'_n\}_{n\in \N}$ based on Lemma \ref{kkk}. In addition, select
$d-2$ integer vectors $l_{n3},\ldots,l_{nd}$ such that
$k_n,k'_n,l_{n3},\ldots,l_{nd}$ are pairwise orthogonal. Let
\begin{equation}\label{KK}
K_n:=(k_n,k'_n,l_{n3},\ldots,l_{nd})^t. \end{equation} We choose the change of the coordinates
\begin{equation}\label{K}
q=K_nx.
\end{equation}
Let $p$ denote the dual coordinate of $q$ in the sense of  Legendre
transformation, i.e. $p=\frac{\partial L}{\partial\dot{q}}$. It
follows that
$y=K^t_np$, where $K^t_n$ denotes the transpose of $K_n$. Let
\[\Phi_n:=\begin{pmatrix}
K_n&\textbf{0}\\
\textbf{0} &K^{-t}_n \end{pmatrix},\]then
$(q,p)^{t}=\Phi_n(x,y)^t$.
Let
\[J_0:=\begin{pmatrix}
\textbf{0}&\textbf{1}\\
-\textbf{1}&\textbf{0}\end{pmatrix},\]
 where $\textbf{1}$ denotes a
$d\times d$ unit matrix.
It is easy to verify that
$\Phi_n^tJ_0\Phi_n=J_0$.
 Hence, $\Phi_n$ is a symplectic
transformation in the phase space $\mathrm{T}^*\mathbb{T}^d$.
 In particular, we have
\begin{equation}\label{cortranfor}
\begin{cases}q_1=\langle k_n,x\rangle,\\
q_2=\langle k'_n,x\rangle.
\end{cases}\end{equation}
In the new coordinates, the rotation vector corresponding to
$\omega$ is given by
$\omega'=K_n\omega$.
 From (\ref{ap}) and (\ref{kn1}), it follows that
\begin{equation}\label{ome1}
|\omega'_1|\lesssim\frac{1}{|k_n|^{d+\tau-1}},\qquad |\omega'_2|\sim |k_n|.
\end{equation}
For the simplicity of notations and without ambiguity, we will still
write the rotation vector by $\omega$ instead of $\omega'$ in the
new coordinates. In addition, we will also use the same notation to denote the rotation vector of the orbit on the torus and its lift to $\R^d$.

\begin{Lemma}\label{key11}
If the Hamiltonian flow generated by $\check{H}_n(x,y)$ admits a
Lagrangian torus with the rotation vector $\omega$, then the Hamiltonian
flow generated by $\hat{H}_n(q,p)$ also admits a Lagrangian torus with the
rotation vector $K_n\omega$, where $(q,p)^t=\Phi_n(x,y)^t$.
\end{Lemma}

\Proof Let $\check{\mathcal {T}}^d$ be the Lagrangian torus admitted
by $\check{H}_n(x,y)$.
Since $K_n$ consists of integer vectors, then $\hat{\mathcal
{T}}^d:=K_n\check{\mathcal {T}}^d$ is still a torus. Note that $\Phi_n$ is a
symplectic transformation,  $\hat{\mathcal {T}}^d$ is a Lagrangian
torus. From Definition \ref{dd1}, there exists a dense orbit $\check{\gm}_0$ with the rotation vector $\omega$ on $\check{\mathcal {T}}^d$. Let
$\check{\gamma}(t)$ be the lift of $\check{\gm}_0$ to $\R^d$, it follows that
\begin{equation}\label{rotv}
\omega=\lim_{t-s\rightarrow+\infty}\frac{\check{\gamma}(t)-\check{\gamma}(s)}{t-s}.
\end{equation}
Let $\hat{\gamma}(t):=K_n\check{\gamma}(t)$, we have
\begin{align*}
\lim_{t-s\rightarrow+\infty}\frac{\hat{\gamma}(t)-\hat{\gamma}(s)}{t-s}&=\lim_{t-s\rightarrow+\infty}\frac{K_n\check{\gamma}(t)-K_n\check{\gamma}(s)}{t-s};\\
&=K_n\lim_{t-s\rightarrow+\infty}\frac{\check{\gamma}(t)-\check{\gamma}(s)}{t-s};\\
&=K_n\omega.
\end{align*}
Let $\hat{\gm}_0$ be the projected orbit of $\hat{\gm}$ on $\hat{\mathcal
{T}}^d$. Note that
\[\overline{\{\hat{\gm}_0(t)\ |\ t\in \R\}}=K_n \overline{\{\check{\gm}_0(t)\ |\ t\in \R\}}.\]
It follows that $\hat{\gm}_0$ is also dense on $\hat{\mathcal
{T}}^d$.
This completes the proof.\End

\section{\sc Construction of the perturbation}\label{3s}

Given $0<\eps\ll 1$, the perturbation  $P_{N}(x)$ is constructed as follows:
\[P_{N}(x)=\frac{1}{|k_n|^{a+2}}(1-\cos \langle k_n,x\rangle)+\frac{1}{|k_n|^{2}}v_n(\langle k_n,x\rangle,\langle k'_n,x\rangle),\] where $k_n$, $k'_n$ are $d$-dimensional vectors given by the first two rows of $K_n$ defined as (\ref{KK}) and $a$ is a parameter given by
\begin{equation}\label{3333}
a:=2d+2\tau-2-\eps.
\end{equation}
 Let $q=K_nx$. We have $q_1=\langle k_n,x\rangle$ and $q_2=\langle k'_n,x\rangle$. Then $v_n(q_1,q_2)$ is a $2\pi$-periodic function with respect to $(q_1,q_2)$.

We construct $v_n(q_1,q_2)$ in the following way. First of all, we choose a  $2\pi$-periodic $C^\infty$ function $\varphi_n$ defined on $\R$. Within a fundamental domain $[-\pi,\pi]$,  we require its maximum $\varphi_n(0)=\sqrt{2}$, and it  is supported on
\[\left[-\frac{|\omega_1|}{|k_n|^{1+\eps}},\frac{|\omega_1|}{|k_n|^{1+\eps}}\right].\]
By (\ref{ome1}), we have
\begin{equation}\label{ommm1}
|\omega_1|\lesssim\left(\frac{1}{|k_n|}\right)^{d+\tau-1}.
\end{equation}
 Let $\phi_n(q_1,q_2):=\varphi_n(q_1-\pi)\cdot\varphi_n(q_2)$. For simplicity, we denote
\[R_n:=\frac{|\omega_1|}{|k_n|^{1+\eps}}.\]Then
\begin{equation}\label{phin}
\left\{\begin{array}{ll}
\hspace{-0.4em}\phi_n(\pi,0)=\max\phi_n(q_1,q_2)=2,\\
\hspace{-0.4em}\text{supp}\phi_n=\left[\pi-R_n,\pi+R_n\right]\times
\left[-R_n,R_n\right],\\
\hspace{-0.4em}\|\phi_n\|_{C^{\kappa}}\sim R_n^{-\kappa}.
\end{array}\right.
\end{equation}
Next, we use Jackson's approximation (\cite{Z}) to obtain a trigonometric polynomial $T_{M,n}(q_1,q_2)$. After the improvement by Favard \cite{Fa}, we have
\begin{Theorem}[Jackson's approximation]Let $\phi(x)$ be a
$\kappa$-times ($\kappa\geq 2$) differentiable $2\pi$-periodic function on $\R$, then for every
$M\in\N$, there exists a trigonometric polynomial $f_M(x)$ of degree
$M$ such that
\[\max|f_M(x)-\phi(x)|\leq \frac{4C_\kappa}{\pi(M+1)^\kappa}||\phi(x)||_{C^{\kappa}},\]
where $C_\kappa=\Sigma_{i=1}^\infty\frac{(-1)^{i(\kappa-1)}}{(2i+1)^{\kappa+1}}$.
\end{Theorem}
Note that $C_\kappa$ is an absolute constant independent of $\phi$ and $f_M$. Hence, for $\varphi_n$ defined as above, there exists a trigonometric polynomial $g_{M,n}(x)$ of degree
$M$ such that
\[\max|g_{M,n}(x)-\varphi_n(x)|\lesssim M^{-\kappa}||\varphi_n(x)||_{C^{\kappa}}.\]
 Let
\[T_{M,n}(q_1,q_2):=g_{M,n}(q_1-\pi)g_{M,n}(q_2).\]
It follows that
\[\max|T_{M,n}(q_1,q_2)-\phi_n(q_1,q_2)|\lesssim M^{-\kappa}||\phi_n||_{C^{\kappa}}.\]
Based on the construction of $\phi_n$, there hold
\begin{equation}\label{Tnn}
\left\{\begin{array}{ll}
\hspace{-0.4em} T_{M,n}(\pi,0)\geq 1,&\\
\hspace{-0.4em} T_{M,n}(q_1,q_2)\leq \mu_n\ll 1, &\forall (q_1,q_2)\in ([0,2\pi]\times[-\pi,\pi])\backslash \text{supp}\phi_n.
\end{array}\right.
\end{equation}
More precisely, we require
\begin{equation}\label{munn}
\mu_n:=M^{-\kappa}||\phi_n||_{C^{\kappa}}\sim |\omega_1|^{\alpha}.
\end{equation}
We choose $\kappa=\frac{\alpha}{\eps}$, where   $\alpha\geq 2$  is a constant independent of $n$ (determined by (\ref{alp}) blew). Combining with (\ref{phin}), we have
\begin{equation}\label{Norde}
M\sim \frac{|k_n|^{1+\eps}}{|\omega_1|^{1-\eps}}.
\end{equation}

By Bernstein's inequality and Leibniz's rule,
for any fixed $s\geq 0$, we have
\begin{equation*}
||T_{M,n}(q_1,q_2)||_{C^s}\lesssim  M^{s}\max|T_{M,n}(q_1,q_2)|.
\end{equation*}
After normalization, we denote
\[\bar{T}_{M,n}(q_1,q_2):=M^{-2s_0}\left(\frac{T_{M,n}(q_1,q_2)}{\max T_{M,n}(q_1,q_2)}\right)^2,\]
where $s_0:=2d+2\tau$.
Finally, we construct $v_n$ as follows.
\begin{equation*}
v_n(q_1,q_2)=\frac{1}{|k_n|^a}(1-\cos q_1)\bar{T}_{M,n}(q_1,q_2).
\end{equation*}
Moreover, there hold
\begin{equation}\label{vnnes}
\left\{\begin{array}{ll}
\hspace{-0.4em} v_n(q_1,q_2)\geq 0,\\
\hspace{-0.4em} \|v_n\|_{C^r}\lesssim M^{-2(s_0-r)}\frac{1}{|k_n|^a}, \quad \forall r< s_0,\\
\hspace{-0.4em} \max v_n\sim M^{-2s_0}\frac{1}{|k_n|^a},\\
\hspace{-0.4em} \|v_n\|_{C^0}\sim \mu_n^2M^{-2s_0}\frac{1}{|k_n|^a}, \quad \forall (q_1,q_2)\in ([0,2\pi]\times[-\pi,\pi])\backslash \text{supp}\phi_n.
\end{array}\right.
\end{equation}

By the construction of $P_{N}(x)$,  for $r< 2d+2\tau-\eps$, we have
\begin{equation}\label{pest}
\|P_{N}\|_{C^r}\sim \left(\frac{1}{|k_n|}\right)^{2d+2\tau-\eps-r}.
\end{equation}
Let ${\delta}:=2\eps$. A direct calculation implies
\begin{equation}\label{cinf}
\|H_n(x,y)-H_0(y)\|_{C^{2d+2\tau-{\delta}}}=\|P_{N}(x)\|_{C^{2d+2\tau-{\delta}}}\sim \left(\frac{1}{|k_n|}\right)^{\eps}\rightarrow 0\quad \text{as}\ n\rightarrow \infty.
\end{equation}
The proof of Item (1) of Theorem \ref{mathe} will be completed, once we  verify Proposition \ref{keyy} below. Here, we first prove Item (2) and (3).

Given $0<\eps\ll 1$, we require $\|P_{N}\|_{C^r}\sim\eps$. It follows from (\ref{pest}) that
\begin{equation}\label{mMM}
|k_n|\sim \eps^{-\frac{1}{2d+2\tau-\eps-r}}.
\end{equation}
Since $\omega$ has the Diophantine exponent $d-1+\tau$, by definition, we have
\[ |\omega_1|\gtrsim\frac{1}{|k_n|^{d-1+\tau+\eps}},\]
which together with (\ref{Norde}) implies
\begin{equation}\label{MM1}
M\lesssim |k_n|^{d+\tau+\eps(3-d-\tau-\eps)}.
\end{equation}
Combining with (\ref{mMM}), we have
\begin{equation}\label{MM}
M\lesssim \eps^{-\frac{d+\tau}{2d+2\tau-\eps-r}-\frac{\eps(3-d-\tau-\eps)}{2d+2\tau-\eps-r}}.
\end{equation}
Since $\eps$ is small enough, (\ref{MM}) can be reduced to
\[M\lesssim\eps^{-\frac{d+\tau}{2d+2\tau-r}}.\]
Note that the degree $N$ of $P_{N}$ is not greater than $(2M+1)|k_n|$. Then  we get
\[N\lesssim\eps^{-\frac{d+\tau+1}{2d+2\tau-r}}.\]

If $\omega$ is   a {Liouvillean type}, by taking $r=\ln \tau$ for instance,  then the Lagrangian torus can be destroyed by  $P_N$ of degree $N\lesssim {{\epsilon}^{-1/2}}$. It remains to verify $\|P_N\|_{C^\infty}\lesssim \epsilon$. In fact, for any $\tau>0$, it follows from (\ref{cinf}) that
\[\|P_{N}\|_{C^{\ln\tau}}\sim \left(\frac{1}{|k_n|}\right)^{2d+2\tau-\eps-\ln\tau}\rightarrow 0\quad \text{as}\ n\rightarrow \infty.\]
We require $\sup_{\tau>0}\|P_{N}\|_{C^{\ln \tau}}\sim\eps$. Then
\[\|P_{N}\|_{C^\infty}=\sum_{k\in \mathbb{N}}\frac{\arctan (\|P_N\|_k)}{2^k}\leq \sup_{\tau>0}\|P_{N}\|_{C^{\ln \tau}}\sim\eps.\]

\section{\sc Lagrangian formulatioin}\label{41s}

According to (\ref{K}), the Lagrangian corresponding to (\ref{Lo}) is
\begin{equation}\label{lnn}
\begin{split}
L_n(q,\dot{q})=&\frac{1}{|k_n|^2}\left(\frac{1}{2}\sum_{i=3}^d\frac{|k_n|^2}{|l_{ni}|^2}|\dot{q}_i|^2+\frac{1}{2}\frac{|k_n|^2}{|k'_n|^2}|\dot{q}_2|^2\right.\\
&\left.+\frac{1}{2}|\dot{q}_1|^2
+\frac{1}{|k_n|^a}(1-\cos q_1)+v_n(q_1,q_2)\right),
\end{split}
\end{equation}
where $a$ is given by (\ref{3333}). Let
\begin{equation}\label{sipp}
A_n(q_1,\dot{q}_1):=\frac{1}{2}|\dot{q}_1|^2
+\frac{1}{|k_n|^a}(1-\cos q_1);
\end{equation}
\begin{equation}\label{bnes}
B_n(q,\dot{q}):=\frac{1}{2}\sum_{i=3}^d\frac{|k_n|^2}{|l_{ni}|^2}|\dot{q}_i|^2+\frac{1}{2}\frac{|k_n|^2}{|k'_n|^2}|\dot{q}_2|^2+v_n(q_1,q_2).
\end{equation}
Then we have
\begin{equation}\label{lnab}
L_n=\frac{1}{|k_n|^2}(A_n+B_n).
\end{equation}

In view of Lemma \ref{key11} and (\ref{ome1}), in order to complete the proof of Theorem \ref{mathe}(1), we only need to prove
\begin{Proposition}\label{keyy}
For $n$ large enough, the Lagrangian system generated by $L_n(q,\dot{q})$ defined as (\ref{lnn}) does not admit the Lagrangian torus with the rotation vector $\omega:=(\omega_1,\omega_2,\dots,\omega_d)$ satisfying
\begin{equation}\label{ap1}
|\omega_1|\lesssim\left(\frac{1}{|k_n|}\right)^{d+\tau-1},\quad {\omega_2}\sim |k_n|.
\end{equation}
\end{Proposition}

As preparations, we recall some facts about the dynamics generated by $A_n$ and $B_n$ respectively. Let $q_1(t)$ be a solution of $A_n$ on $(t_0,t_1)$ and $(t_1,t_2)$
with boundary conditions respectively
\begin{equation*}\begin{cases}q_1(t_0)=0,\\
q_1(t_1)=\pi,\end{cases}
\begin{cases}q_1(t_1)=\pi,\\
q_1(t_2)=2\pi.
\end{cases}\end{equation*} We consider the function $\mathbb{L}:[t_0,t_2]\to\R$ defined by
\begin{align*}
\mathbb{L}(s):=&\int_{t_0}^{s}A_n(q_1,\dot{q}_1)dt+\int_{s}^{t_2}A_n(q_1,\dot{q}_1)dt.
\end{align*}
By \cite[Remark 2.5]{CW}, there hold
\begin{Lemma}\label{AR}
$\mathbb{L}(s)$ is decreasing for $s\in
(t_0,\frac{t_0+t_2}{2}]$ and increasing for $s\in
[\frac{t_0+t_2}{2},t_2)$.
\end{Lemma}

Another observation is that the actions along orbits in the
neighborhood of the separatix of the  pendulum do not change too
much with respect to a small change in speed. To set up a precise statement, we introduce some notations.
Given
$\bar{t}_1, \tilde{t}_1\in [t_0,t_2]$. Let $\bar{q}_1(t)$ be a
solution of $A_n$ on $(t_0,\bar{t}_1)$ and $(\bar{t}_1,t_2)$ with
boundary conditions respectively
 \begin{equation*}\begin{cases}\bar{q}_1(t_0)=0,\\
\bar{q}_1(\bar{t}_1)=\pi,\end{cases}
\begin{cases}\bar{q}_1(\bar{t}_1)=\pi,\\
\bar{q}_1(t_2)=2\pi.\end{cases}\
\end{equation*}
Let $\tilde{q}_1(t)$ be a solution of $A_n$ on
$(t_0,\tilde{t}_1)$ and $(\tilde{t}_1,t_2)$ with boundary conditions
respectively
\begin{equation*}
\begin{cases}\tilde{q}_1(t_0)=0,\\
\tilde{q}_1(\tilde{t}_1)=\pi,\end{cases}
\begin{cases}\tilde{q}_1(\tilde{t}_1)=\pi,\\
\tilde{q}_1(t_2)=2\pi.
\end{cases}
\end{equation*} Let
$\bar{\omega}_1$ and $\bar{\omega}_2$ be the average speeds of
$\bar{q}_1$ on $(t_0,\bar{t}_1)$ and $(\bar{t}_1,t_2)$ respectively.
Let $\tilde{\omega}_1$ and $\tilde{\omega}_2$ be the average
speed of $\tilde{q}_1$ on $(t_0,\tilde{t}_1)$ and
$(\tilde{t}_1,t_2)$ respectively. By \cite[Lemma 2.4]{CW}, there hold
\begin{Lemma}\label{8_1}
Let
\[\omega_1:=\max\left\{\bar{\omega}_1,\bar{\omega}_2,\tilde{\omega}_1,\tilde{\omega}_2\right\}.\] If
\[0<\omega_1\lesssim\left(\frac{1}{|k_n|}\right)^{d+\tau-1},\]
then
\begin{equation}
\left|\int^{t_2}_{t_0}A_n(\bar{q}_1,\dot{\bar{q}}_1)dt-
\int^{t_2}_{t_0}A_n(\tilde{q}_1,\dot{\tilde{q}}_1)dt\right|\lesssim\frac{|\bar{t}_1-\tilde{t}_1|}{|k_n|^{a}}
\exp\left(-\frac{C}{\omega_1|k_n|^{\frac{a}{2}}}\right),
\end{equation}
where $a$  is given by (\ref{3333}).
\end{Lemma}

The last fact is about the speed of the action minimizing orbit of $L_n$. Once the function $q_1(t)$ is fixed, we denote $Q(t):=(q_2(t),\ldots,q_d(t))$ to be the
solution of the Euler-Lagrange equation with the non autonomous
Lagrangian
\begin{equation}\label{LQ}
\hat{B}_n(Q(t),\dot{Q}(t),t)=\frac{1}{2}\sum_{i=3}^d\frac{|k_n|^2}{|l_{ni}|^2}|\dot{q}_i|^2+\frac{1}{2}\frac{|k_n|^2}{|k'_n|^2}|\dot{q}_2|^2+v_n\left(q_1(t),
q_2(t)\right),
\end{equation}The first order
derivative of $\hat{B}_n(Q(t),\dot{Q}(t),t)$ with respect to $Q$ is as follows:
\begin{equation}\label{eleq}
\frac{\partial \hat{B}_n}{\partial q_2}=\frac{\partial v_n}{\partial q_2}(q_1(t),q_2(t))\quad\text{and}
\quad\frac{\partial \hat{B}_n}{\partial q_i}=0\quad\text{for}\quad
i=3,\ldots,d.
\end{equation} From the construction of $v_n$ (see (\ref{vnnes})), we
have $||v_n||_{C^r}\lesssim\frac{1}{|k_n|^a}$ if $r< 2d+2\tau$. Recalling $a:=2d+2\tau-2-\eps$, it follows from (\ref{eleq}) that
\[\left|\frac{\partial \hat{B}_n}{\partial q_i}\right|\lesssim\frac{1}{|k_n|^a}\quad\text{for}\quad
i=2,\ldots,d..\] Note that
$\hat{B}_n(Q(t),\dot{Q}(t),t)$ is periodic with respect to $Q(t)$, and $a>1-\frac{\eps}{2}$ for all $d\geq 2$ and $\tau\geq 0$. By   \cite[Lemma 2]{BK}, we have the
following estimate.

\begin{Lemma}\label{9_1} Let $(q_1(t),Q(t))$ be the  action minimizing orbit of $L_n$ with the
rotation vector $\omega$ satisfying (\ref{ap1}), then for any $t',t''\in \R$ and $t\in
[t',t'']$ we have
\begin{equation}\label{vv}
\left|\dot{Q}(t)-\frac{Q(t'')-Q(t')}{t''-t'}\right|\lesssim\frac{1}{|k_n|^{d+\tau-1-\frac{\eps}{2}}}\lesssim\frac{1}{\sqrt[3]{|k_n|}}.
\end{equation}\end{Lemma}

\section{\sc Proof of Proposition \ref{keyy}}\label{4s}
It is sufficient to prove the existence of the small neighborhood of some point in
 $\mathbb{T}^d:=\R^d/ 2\pi\mathbb{Z}^d$ where no action minimizing curve passes through.
Let us recall
\[\text{supp}\phi_n=\left[\pi-R_n,\pi+R_n\right]\times
\left[-R_n,R_n\right],\]
where
\[R_n:=\frac{|\omega_1|}{|k_n|^{1+\eps}}.\]
We denote
\[S_0:=\left[\pi-\frac{1}{2}R_n,\pi+\frac{1}{2}R_n\right]\times
\left[-\frac{1}{2}R_n,\frac{1}{2}R_n\right]\]
 In fact, we will show that the action minimizing orbit does not pass through
\[S_0\times \T^{d-2}.\]
By contradiction, we assume that there exists $\bar{t}_1$ such that
\begin{equation}\label{qqx}
(q_1(\bar{t}_1),q_2(\bar{t}_1))\in S_0.
\end{equation} Without loss of generality, we assume
\[q_1(\bar{t}_1)=\pi,\quad q_2(\bar{t}_1)=0,\quad \omega_1>0,\quad \omega_2>0,\] where
$q(t):=(q_1,q_2,\ldots,q_d)(t)$ is an action minimizing orbit in the universal
covering space $\R^d$. Moreover, there exist $t_0<t_2$ such that
\[q_1(t_0)=0,\quad q_1(t_2)=2\pi.\]

\vspace{1ex}

\noindent\textbf{Claim:}
$t_2-t_0\sim {{\omega_1}}^{-1}$.

\noindent\Proof By contradiction, we assume $t_2-t_0\sim {\bar{\omega}_1}^{-1}$, and (up to a subsequence) either ${\bar{\omega}_1}=o({\omega_1})$, or ${\omega_1}=o({\bar{\omega}_1})$. We only need to consider the case with ${\bar{\omega}_1}=o({\omega_1})$, since the other case can be excluded by a similar argument.
The aim is to construct another curve with smaller action compared to the orbit $q(t)$. That contradicts the action minimizing property of $q(t)$.

By the definition of the rotation vector,  for any $\eps>0$, there exists $t_N>t_2$ such that
\begin{equation}
\left|\frac{q_1(t_N)-q_1(\bar{t}_1)}{t_N-\bar{t}_1}-\omega_1\right|\leq\epsilon\quad\text{and}
\quad q_1(t_N)=N\pi,
\end{equation}where $N$ depends on $n$ and $N\gg 2$. Then we have
\[t_N-\bar{t}_1\sim \frac{N\pi}{\omega_1}.\]We choose a sequence of times $t_2<\ldots<t_N$ satisfying $q_1(t_i)=i\pi$ for $i\in
\{2,\ldots,N\}$. Moreover, it follows from the pigeon hole principle
that there exists $j\in \{2,\ldots,N-1\}$ such that
\[t_{j+1}-t_j\lesssim\frac{1}{\omega_1},\]where we consider the case
$q_1(t_j)\text{mod}2\pi=0$ and $q_1(t_{j+1})\text{mod}2\pi=\pi$, the
other case is similar. Let $q_1(\bar{t}'_1)=\pi+2R_n$ and
$q_1(t'_{j+1})=(j+1)\pi-2R_n$.

The remaining proof of this claim is divided into three steps. First of all,  we show $t_2-\bar{t}'_1\sim {\bar{\omega}_1}^{-1}$. Denote the average speeds of $q_1(t)$ passing through $[\pi, \pi+2R_n]$ and $[\pi+2R_n,2\pi]$ by $\omega'_1$ and $\omega''_1$ respectively. Note that
\begin{equation}\label{prot}
\frac{\pi}{\bar{\omega}_1}=\frac{2R_n}{\omega'_1}+\frac{\pi-2R_n}{\omega''_1}.
\end{equation}
It follows that
\[\omega''_1\geq \frac{\pi-2R_n}{\pi}\bar{\omega}_1.\]
By the Euler-Lagrange equation generated by $L_n$, it is direct to see $\omega'_1>\omega''_1$. Then
\[\bar{\omega}_1\lesssim \omega''_1\leq \bar{\omega}_1,\]
which implies $t_2-\bar{t}'_1\sim {\bar{\omega}_1}^{-1}$.

Second, we  substitute
$q_1(t)|_{[\bar{t}'_1,t_2]}$ and $q_1(t)|_{[t_j,t'_{j+1}]}$ by the
orbits $\hat{q}_{1'}(t)$, $\hat{q}_{1''}(t)$ of the pendulum $A_n$
(see (\ref{sipp})) with more uniform motion and smaller action (see Fig. 1).
\begin{figure}[htbp]\label{fi3}
\small \centering
\includegraphics[width=14.0cm]{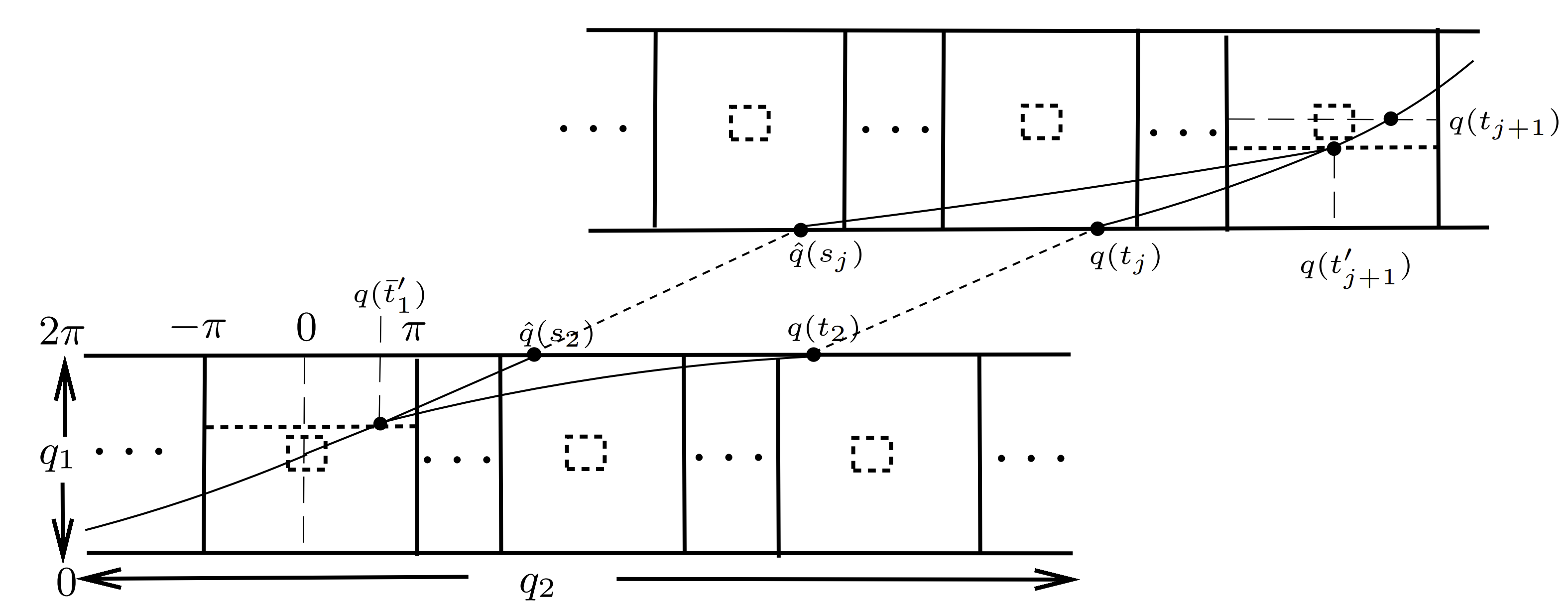}
\caption{Construction of the curve with smaller action}
\end{figure}
In fact, by a similar argument as the first step, we have $t'_{j+1}-t_j\sim {\omega_1}^{-1}$. Since we assume  ${\bar{\omega}_1}=o({\omega_1})$,  the second step can be achieved by using Lemma \ref{AR}. More precisely,
we denote
\[T=\frac{1}{2}(t_2-\bar{t}'_1+t'_{j+1}-t_j).\]
Let $s_{j}$ be the closest time to $t'_{j+1}-T$ such that
$q_2(s_j)-q_2(t_j)= l$ where $l\in 2\pi\Z$. Since ${\omega_2}\sim
|k_n|$, it follows from Lemma \ref{9_1} that $|t'_{j+1}-T-s_j|\lesssim |k_n|^{-1}$.
Let $s_2=s_j-(t_j-t_2)$. Then $\hat{q}_{1'}(t)$ and
$\hat{q}_{1''}(t)$ are the solutions of $A_n$ on $(\bar{t}'_1,s_2)$
and on $(s_j,t'_{j+1})$ with boundary conditions respectively
\begin{equation*}
\begin{cases}
\hat{q}_{1'}(\bar{t}'_1)=q_1(\bar{t}'_1),\\
\hat{q}_{1'}(s_2)=q_1(t_2)=2\pi,
\end{cases}
\quad
\begin{cases}
\hat{q}_{1''}(s_j)=q_1(t_j)=j\pi,\\
\hat{q}_{1''}(t'_{j+1})=q_1(t'_{j+1}).
\end{cases}
\end{equation*}

Finally, we substitute
$q_2(t)|_{[\bar{t}'_1,t_2]}$ and $q_2(t)|_{[t_j,t'_{j+1}]}$ by the linear motions $\hat{q}_{2'}(t)$ and
$\hat{q}_{2''}(t)$ respectively. On the time interval $[\bar{t}'_1, t'_{j+1}]$, we construct a new curve as follows.
\[\hat{q}_1:=\hat{q}_{1'}(t)|_{[\bar{t}'_1,s_2]}\ast q_1(t-s_2+t_2)|_{[s_2,s_j]}\ast
\hat{q}_{1''}(t)|_{[s_j,t'_{j+1}]},\]
\[\hat{q}_2:=\hat{q}_{2'}(t)|_{[\bar{t}'_1,s_2]}\ast \left(q_2(t-s_2+t_2)|_{[s_2,s_j]}-l\right)
\ast \hat{q}_{2''}(t)|_{[s_j,t'_{j+1}]},\] where $\ast$ denotes the
juxtaposition of curves. To complete the proof of this claim, it suffices to verify
\begin{equation}
\int_{\bar{t}'_1}^{ t'_{j+1}}L_n(\hat{q}_1,\hat{Q},\dot{\hat{q}}_1,\dot{\hat{Q}})dt<\int_{\bar{t}'_1}^{ t'_{j+1}}L_n(q_1,Q,\dot{q}_1,\dot{Q})dt,
\end{equation}
where $\hat{Q}:=(\hat{q}_2,q_3,\ldots,q_d)$
and $Q:=(q_2,q_3,\ldots,q_d)$. Let us recall (see (\ref{lnab}))
\[L_n=\frac{1}{|k_n|^2}(A_n+B_n).\]
Based on the second step, we only need to show
\begin{equation}\label{aibb}
\begin{split}
&\int_{\bar{t}'_1}^{t_2}B_n(\hat{q}_1,\hat{Q},\dot{\hat{q}}_1,\dot{\hat{Q}})dt+\int_{t_j}^{t'_{j+1}}B_n(\hat{q}_1,\hat{Q},\dot{\hat{q}}_1,\dot{\hat{Q}})dt\\
\leq &\int_{\bar{t}'_1}^{t_2}B_n(q_1,Q,\dot{q}_1,\dot{Q})dt+\int_{t_j}^{t'_{j+1}}B_n({q}_1,{Q},\dot{{q}}_1,\dot{{Q}})dt.
\end{split}
\end{equation}
Denote
\[d_1:=q_2(t_2)-q_2(\bar{t}'_1),\quad d_2:=q_2(t'_{j+1})-q_2(t_j).\]
According to the construction of $v_n$, both actions of $q(t)$ and $\hat{q}(t)$ could be approximated by the actions of linear motions. It follows that
\[\int_{\bar{t}'_1}^{t_2}B_n(q_1,Q,\dot{q}_1,\dot{Q})dt+\int_{t_j}^{t'_{j+1}}B_n({q}_1,{Q},\dot{{q}}_1,\dot{{Q}})dt\sim \frac{d_1^2}{t_2-\bar{t}'_1}+\frac{d_2^2}{t'_{j+1}-t_j}.\]
\[\int_{\bar{t}'_1}^{t_2}B_n(\hat{q}_1,\hat{Q},\dot{\hat{q}}_1,\dot{\hat{Q}})dt+\int_{t_j}^{t'_{j+1}}B_n(\hat{q}_1,\hat{Q},\dot{\hat{q}}_1,\dot{\hat{Q}})dt\sim \frac{(d_1-l)^2}{s_2-\bar{t}'_1}+\frac{(d_2+l)^2}{t'_{j+1}-s_j}.\]
Note that $|t'_{j+1}-T-s_j|\lesssim |k_n|^{-1}$. The average speeds of $\hat{q}_1'$ and $\hat{q}_1''$ have the same quantity order. Let $\hat{\omega}_1$ be the average speed of $\hat{q}_1'$ and $\hat{q}_1''$. Similar to (\ref{prot}), we have
\[\hat{\omega}_1\sim \min\{\omega_1,\bar{\omega}_1\}=\bar{\omega}_1,\]
which yields
\[s_2-\bar{t}'_1\sim t'_{j+1}-s_j\sim \frac{1}{\hat{\omega}_1}.\]
Combining with the following facts
\[t_2-\bar{t}'_1\sim \frac{1}{\bar{\omega}_1},\quad t'_{j+1}-t_j\lesssim\frac{1}{{\omega}_1},\quad d_1-l,d_2+l\in [d_2,d_1],\quad \bar{\omega}_1=o(\omega_1),\]
We have
\begin{equation}
\begin{split}
\frac{(d_1-l)^2}{s_2-\bar{t}'_1}+\frac{(d_2+l)^2}{t'_{j+1}-s_j}&\lesssim \bar{\omega}_1(d_1^2+d_2^2)=o\left(\bar{\omega}_1d_1^2+\omega_2d_2^2\right),\\
\frac{d_1^2}{t_2-\bar{t}'_1}+\frac{d_2^2}{t'_{j+1}-t_j}&\succeq\bar{\omega}_1d_1^2+\omega_2d_2^2.
\end{split}
\end{equation}
This completes the proof of the claim $t_2-t_0\sim {{\omega_1}}^{-1}$.
\End

 Let $\tilde{t}_1$ be the last time before $\bar{t}_1$ or the first
time after $\bar{t}_1$ such that
\[|q_2(\tilde{t}_1)-q_2(\bar{t}_1)|=\pi.\]
 Consider a solution $\tilde{q}_1$ of
$A_n$ on $(t_0,\tilde{t}_1)$ and on $(\tilde{t}_1,t_2)$ with
boundary conditions respectively
\begin{equation*}
\begin{cases}
\tilde{q}_1(t_0)=0,\\
\tilde{q}_1(\tilde{t}_1)=\pi,
\end{cases}
\quad
\begin{cases}
\tilde{q}_1(\tilde{t}_1)=\pi,\\
\tilde{q}_1(t_2)=2\pi.
\end{cases}
\end{equation*}
Note that $q=(q_1,q_2,\ldots,q_d)$ is assumed to be an action minimizing curve.
Let $Q:=(q_2,\ldots,q_d)$, we have
\begin{equation}\label{ssg}\int_{t_0}^{t_2}L_n(\tilde{q}_1,Q,\dot{\tilde{q}}_1,\dot{Q})dt-\int_{t_0}^{t_2}L_n(q_1,Q,\dot{q}_1,\dot{Q})dt\geq
0.\end{equation}
See Fig.2, where
$x_1=(q_1(\bar{t}_1),q_2(\bar{t}_1))=(\pi,0)$, $x_0=(q_1(t_0),q_2(t_0))=(0,q_2(t_0))$, $x_2=(q_1(t_2),q_2(t_2))=(2\pi,q_2(t_2))$,
$\tilde{x}'_1=(\pi,-\pi)$,  $\tilde{x}''_1=(\pi,\pi)$.

\begin{figure}[htbp]\label{fi2}
\small \centering
\includegraphics[width=13cm]{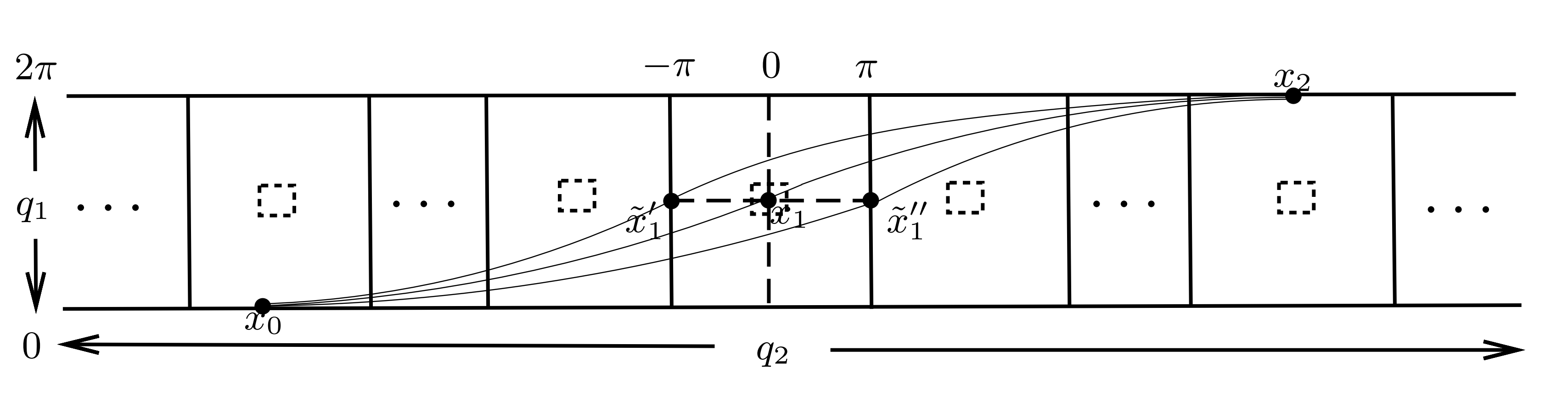}
\caption{The projections of the curves $(q_1(t),Q(t))$ and
$(\tilde{q}_1(t),Q(t))$ on $[0,2\pi]\times\R$}
\end{figure}

Note that $(\tilde{q}_1(t),q_2(t))$ passes through the point $\tilde{x}'_1$ or
$\tilde{x}''_1$.  Thus, by the construction of $L_n$, it follows from
(\ref{ssg}) that
\begin{equation}\label{AandP_1ge}
\int_{t_0}^{t_2}A_n(\tilde{q}_1,\dot{\tilde{q}}_1)dt-\int_{t_0}^{t_2}A_n(q_1,\dot{q_1})dt\geq
\int_{t_0}^{t_2}v_n(q_1,q_2)dt-\int_{t_0}^{t_2}v_n(\tilde{q}_1,q_2)dt.
\end{equation}
Based on the construction of $v_n$ (see (\ref{vnnes})), we have the following claim.

\vspace{1ex}

\noindent\textbf{Claim:}
$\left(\tilde{q}_1(t),q_2(t)\right)\notin\text{supp}\phi_n$ for any $t\in (t_0,t_2)$.

\noindent\Proof By contradiction, we assume that there exists $\hat{t}\in (t_0,t_2)$ such that
\begin{equation}\label{inphi}
(\tilde{q}_1(\hat{t}),q_2(\hat{t}))\in
\text{supp}\phi_n.
\end{equation} Without loss of generality, one can assume $\hat{t}>\tilde{t}_1$.  Note that ${\omega_2}\sim |k_n|$. By Lemma \ref{9_1}, We obtain  that for $t\in
[\tilde{t}_1,\hat{t}]$,
\[\dot{q}_2(t)\sim |k_n|.\]
 Hence, we have
$\hat{t}-\tilde{t}_1\sim |k_n|^{-1}$.

Let $\tilde{\omega}'_1$ and $\tilde{\omega}''_1$ be the average
speeds of $\tilde{q}_1$ on $(t_0,\tilde{t}_1)$ and
$(\tilde{t}_1,t_2)$ respectively, then
\[\frac{2\pi}{{\omega_1}}=\frac{\pi}{\tilde{\omega}'_1}+\frac{\pi}{\tilde{\omega}''_1}.\]Hence,
\[\tilde{\omega}'_1\geq \frac{1}{2}{\omega_1}\quad\text{and}\quad\tilde{\omega}''_1\geq \frac{1}{2}{\omega_1},\]
which together with the Euler-Lagrange equation of $\tilde{q}_1(t)$
implies that for any $t\in [\tilde{t}_1,\hat{t}]$,
\[\dot{\tilde{q}}_1(t)\gtrsim{\omega_1},\]
Consequently, we have
\[|\tilde{q}_1(\hat{t})-\tilde{q}_1(\tilde{t}_1)|\gtrsim\frac{{\omega_1}}{|k_n|}.\]It follows that  for $n$ large enough,
\[|\tilde{q}_1(\hat{t})-\tilde{q}_1(\tilde{t}_1)|\gg\frac{{\omega_1}}{|k_n|^{1+\eps}}.\]
 Based on the construction of $\phi_n$ (see (\ref{phin})), that contradicts  the assumption (\ref{inphi}).\End

\vspace{1ex}

\noindent\textbf{Claim:}
There exists $\lambda>0$ independent of $n$ such that
\[\int_{t_0}^{t_2}v_n(q_1,q_2)dt-\int_{t_0}^{t_2}v_n(\tilde{q}_1,q_2)dt\gtrsim\omega_1^\lambda.\]

\noindent\Proof
Based on (\ref{vnnes}) and $t_2-t_0\sim {{\omega_1}^{-1}}$, we need to estimate the lower bound of $\int_{t_0}^{t_2}v_n(q_1,q_2)dt$ and the upper bound of $\int_{t_0}^{t_2}v_n(\tilde{q}_1,q_2)dt$.

By Lemma \ref{9_1}, We obtain  that for $t\in
[{t}_0,{t}_2]$,
\[\dot{q}_2(t)\sim |k_n|.\]
Denote $\triangle t$ the time during which the action minimizing orbit $q(t)$ passes though
\[S_0:=\left[\pi-\frac{1}{2}R_n,\pi+\frac{1}{2}R_n\right]\times
\left[-\frac{1}{2}R_n,\frac{1}{2}R_n\right].\]
It follows that
\[\triangle t\sim \frac{R_n}{|k_n|}.\]
Note that $v_n\geq 0$. By the mean value theorem, we have
\[\int_{t_0}^{t_2}v_n(q_1,q_2)dt\gtrsim M^{-2s_0}\frac{1}{|k_n|^a}\triangle t\sim M^{-2s_0}\frac{{\omega_1}}{|k_n|^{a+2+\eps}},\]
where $s_0:=2d+2\tau$.
Denote ${v}_n^*$ the maximum of $v_n$ restricted on $([0,2\pi]\times[-\pi,\pi])\backslash \text{supp}\phi_n$. Then
\[\int_{t_0}^{t_2}v_n(\tilde{q}_1,q_2)dt\lesssim \frac{{v}_n^*}{{\omega_1}}\sim M^{-2s_0}\frac{1}{|k_n|^a}\omega_1^{2\alpha-1}.\]
It follows that
\[\int_{t_0}^{t_2}v_n(q_1,q_2)dt-\int_{t_0}^{t_2}v_n(\tilde{q}_1,q_2)dt\gtrsim M^{-2s_0}\frac{\omega_1}{|k_n|^a}\left(\frac{1}{|k_n|^{2+\eps}}-\omega_1^{2\alpha-2}\right).\]
Note that (see (\ref{Norde}) and (\ref{ap1}))
\[{\omega_1}\lesssim\left(\frac{1}{|k_n|}\right)^{d+\tau-1},\quad M\sim \frac{|k_n|^{1+\eps}}{\omega_1^{1-\eps}}.\]
We choose
\begin{equation}\label{alp}
\alpha\geq 2,
\end{equation} for which  the claim holds true with $\lambda\geq 2(s_0+6)=4d+4\tau+8$.
\End

 It follows from (\ref{AandP_1ge}) that
\begin{equation}\label{contr1_1ge}
\int_{t_0}^{t_2}A_n(\tilde{q}_1,\dot{\tilde{q}}_1)dt-\int_{t_0}^{t_2}A_n(q_1,\dot{q}_1)dt\succeq
\omega_1^\lambda.
\end{equation}

On the other hand, we consider the solution  $\bar{q}_1(t)$  of $A_n$ on $(t_0,\bar{t}_1)$  and $(\bar{t}_1,t_2)$ with the boundary conditions respectively
\begin{equation*}
\begin{cases}
\bar{q}_1(t_0)=0,\\
\bar{q}_1(\bar{t}_1)=\pi,
\end{cases}
\quad
\begin{cases}
\bar{q}_1(\bar{t}_1)=\pi,\\
\bar{q}_1(t_2)=2\pi,
\end{cases}
\end{equation*}
where $\bar{t}_1$ is determined by (\ref{qqx}). The remaining argument is similar to \cite{CW} with a slight modification. We repeat it here for the reader's convenience.
For $t\in (t_0, \bar{t}_1)$ and $(\bar{t}_1,t_2)$ respectively, the
action of $A_n$ achieves the  minimum along $\bar{q}_1(t)$. Thus, we
have
\[\int_{t_0}^{t_2}A_n(q_1,\dot{q}_1)dt\geq\int_{t_0}^{t_2}A_n(\bar{q}_1,\dot{\bar{q}}_1)dt.\]
Based on the choices of $\tilde{t}_1$, we compare the action
$\int_{t_0}^{t_2}A_n(\tilde{q}_1,\dot{\tilde{q}}_1)dt$ with the
action $\int_{t_0}^{t_2}A_n(\bar{q}_1,\dot{\bar{q}}_1)dt$  in the
following dichotomy. Let us recall  $\bar{t}:=\frac{t_0+t_2}{2}$.

\vspace{5pt}

\noindent Case 1: $|\bar{t}_1-\bar{t}|\lesssim\frac{1}{|k_n|}$.

In this case, the average speeds of $\bar{q}_1$ on $(t_0,\bar{t}_1)$
and $(\bar{t}_1,t_2)$ have the same quantity order as ${\omega_1}$.
By $|\tilde{t}_1-\bar{t}_1|\lesssim \frac{1}{|k_n|}$, we have
$|\tilde{t}_1-\bar{t}|\lesssim \frac{1}{|k_n|}$. Hence the average speeds
of $\tilde{q}_1$ on $(t_0,\tilde{t}_1)$ and $(\tilde{t}_1,t_2)$ have
also the same quantity order as ${\omega_1}$. Thus, Lemma \ref{8_1}
implies
\begin{equation}
\begin{split}
\int_{t_0}^{t_2}A_n(\tilde{q}_1,\dot{\tilde{q}}_1)dt-\int_{t_0}^{t_2}A_n(\bar{q}_1,\dot{\bar{q}}_1)dt&\lesssim
\frac{1}{|k_n|^{a+1}}
\exp\left(-\frac{C}{\omega_1|k_n|^{\frac{a}{2}}}\right).
\end{split}
\end{equation}

\noindent Case 2: $|\bar{t}_1-\bar{t}|>\frac{C}{|k_n|}$ for any positive constant $C$ independent of $n$.

In this case, we take $\tilde{t}_1$ such that
$|\tilde{t}_1-\bar{t}|\leq |\bar{t}_1-\bar{t}|$, which can be
achieved by the suitable choice of the position of
$\tilde{q}_1(\tilde{t}_1)$. More precisely,
\begin{itemize}
  \item [(2a)] if $\bar{t}_1>\bar{t}+\frac{C}{|k_n|}$, we choose $\tilde{t}_1$ as the
last time before $\bar{t}_1$, corresponding to
$(\tilde{q}_1(\tilde{t}_1), q_2(\tilde{t}_1))=\tilde{x}'_1$ in
Fig.2;
  \item [(2b)] if
$\bar{t}_1<\bar{t}-\frac{C}{|k_n|}$, we choose
$\tilde{t}_1$ as the first time after $\bar{t}_1$, corresponding to
$(\tilde{q}_1(\tilde{t}_1), q_2(\tilde{t}_1))=\tilde{x}''_1$ in
Fig.2.
\end{itemize}
 For Case 2a, $\tilde{t}_1\in [\bar{t},\bar{t}_1]$. For Case 2b,  $\tilde{t}_1\in [\bar{t}_1,\bar{t}]$. From Lemma \ref{AR}, it
 follows that $\mathbb{L}(\tilde{t}_1)-\mathbb{L}(\bar{t}_1)\leq 0$, i.e.
\[\int_{t_0}^{t_2}A_n(\tilde{q}_1,\dot{\tilde{q}}_1)dt-\int_{t_0}^{t_2}A_n(\bar{q}_1,\dot{\bar{q}}_1)dt\leq
0.\] Consequently, for any $\bar{t}_1\in (t_0,t_2)$, we can always find
$\tilde{t}_1\in (t_0,t_2)$ such that
\begin{align*}
\int_{t_0}^{t_2}A_n(\tilde{q}_1,\dot{\tilde{q}}_1)dt-\int_{t_0}^{t_2}A_n(q_1,\dot{q}_1)dt&\leq
\int_{t_0}^{t_2}A_n(\tilde{q}_1,\dot{\tilde{q}}_1)dt-\int_{t_0}^{t_2}A_n(\bar{q}_1,\dot{\bar{q}}_1)dt,\\
&\lesssim
\frac{1}{|k_n|^{a+1}}\exp\left(-\frac{C}{\omega_1|k_n|^{\frac{a}{2}}}\right),
\end{align*}
where $a=2d+2\tau-2-\eps$. It contradicts (\ref{contr1_1ge}) for $n$ large enough. Then we completes the proof of Proposition \ref{keyy}. \End

\section{\sc Proof of Theorem \ref{excc}}\label{66}
Let $\rho$ be a Diophantine type rotation vector with the exponent $d-1+\tau$. Let $\mathcal{V}$ be  the set of $C^\omega$ potentials $V:\mathbb{T}^d\to\R$. For each $V\in \mathcal{V}$, we denote $\Phi_V^t$ the flow generated by $H_0+V$. Fixing $r\geq 0$, we denote
\[\mathcal{U}^r_\delta:=\{V\in \mathcal{V}\ |\ \|V\|_{C^r}<\delta\}.\]
It is clear to see that  $\mathcal{U}^r_\delta$ is an open and connected neighborhood of $V\equiv 0$ in the $C^r$ topology.  Following the KAM method developed by Salamon and Zehnder \cite{SaZ} (see also \cite{Sa}), we denote the set of the embedding by
\[\mathcal{W}^k:=\{w=(u,v): \mathbb{T}^d\to \mathbb{T}^d\times\R^d\},\]
where $u$ is a $C^k$ diffeomorphism of $\T^d$, and the embedded torus $w(\T^d)$
 is a $C^k$ graph in $\mathrm{T}^\ast\T^d\cong\T^d\times\R^d$ formed by $\{(x,v\circ u^{-1}(x))\ |\ x\in \T^d\}$. It is known that $w(\T^d)$ consists of quasiperiodic solutions of the Hamiltonian system generated by $H_0+V$. Let
 \[\mathcal{U}(\rho,r,\delta):=\{V\in \mathcal{U}^r_\delta\ |\ \exists \ w\in \mathcal{W}^\omega\ \ s.t. \ \ \Phi_V^t=w\circ R_\rho^t\circ w^{-1}\},\]
 where $R_\rho^t(x)=x+\rho t $ (mod $\mathbb{Z}^d$) is the linear flow on $\T^d$.
We assume (once and for all) $\rho$ is a Diophantine type rotation vector with exponent $d-1+\tau$. According to \cite[Theorem, p89]{SaZ}, there holds
\begin{Proposition}\label{vww}
There exists $\delta>0$ such that for all $s>2d+2\tau$, $\mathcal{U}(\rho,s,\delta)$ is an open set in the $C^\omega$ topology.
\end{Proposition}
Let us consider
\[W_\rho:=\cup_{s>2d+2\tau}\mathcal{U}(\rho,s,\delta).\]
By Proposition \ref{vww}, $W_\rho$ is also an open set of  $\mathcal{U}^r_\delta$ in the $C^\omega$ topology.
Next, we consider the closure of $W_\rho$ (denoted by $\overline{W_\rho}$) in $\mathcal{U}^r_\delta$ with respect to the $C^\omega$ topology.  By a compactness argument, we have
\begin{Proposition}\label{aacoli}
For each $V\in \overline{W_\rho}$, there exists a unique $w\in \mathcal{W}^{\mathrm{lip}}$ such that
\[\Phi_{V}^t=w\circ R_\rho^t\circ w^{-1}.\]
\end{Proposition}
\Proof
 For each $V\in \overline{W_\rho}$, there exist $\{V_n\}_{n\in \mathbb{N}}$ with $V_n\in W_\rho$ such that $V_n$ converges uniformly to $V$ in the $C^\omega$ topology.
By Proposition \ref{vww}, there exists $w_n=(u_n,v_n)\in \mathcal{W}^\omega$ such that
\[\Phi_{V_n}^t=w_n\circ R_\rho^t\circ w_n^{-1}.\]
Regarding $u_n$, the KAM theorem  implies $u_n$ is close to identity up to a translation on $\T^d$ (see \cite[Theorem 2]{Sa}). More precisely,  there exists a constant $C>0$ independent of $n$ such that
\begin{equation}
\|u_n-Id\|_{C^1}\leq C.
\end{equation}
By the Arzel\`{a}-Ascoli theorem,  $u_n$ converges uniformly to a Lipschitz map $u$  (up to a subsequence).

On the other hand, we show the compactness of the family of $\{v_n\}_{n\in\mathbb{N}}$. By the definition of $v_n$, the KAM torus $w_n(\T^d)$ with rotation vector $\rho$ generated by $H_0+V_n$ is formed by $\{(x,v_n\circ u_n^{-1}(x))\ |\ x\in \T^d\}$. Since $H_0+V_n$ is a Tonelli Hamiltonian, a celebrated result in the Aubry-Mather theory states that  the Lipschitz constant of $v_n\circ u_n^{-1}$ only dependes on the $C^2$-norm of the Hamiltonian (\cite[Theorem 2]{M5}). It follows from the assumption that the family of $\{v_n\}_{n\in\mathbb{N}}$ is equi-Lipschitz. Let $(x_n(\cdot),p_n(\cdot)):\R\to \T^d\times\R^d$ be an orbit on $w_n(\T^d)$. Then for all $t\in \R$,
 \[p_n(t)=v_n\circ u_n^{-1}(x_n(t)).\]
 Since $\rho$ is the rotation vector of $w_n(\T^d)$, then there exists $t_0\in \R$ such that \[p_n(t_0)=\dot{x}_n(t_0)=\rho.\] Due to the energy conservation, we have for all $t\in \R$,
\[\frac{1}{2}|p_n(t)|^2+V_n(x_n(t))=\frac{1}{2}|\rho|^2+V_n(x_n(t_0)),\]
which gives rise to
\[\|v_n\circ u_n^{-1}\|_{C^0}\leq\sqrt{\rho^2+4 \|V_n\|_{C^0}}.\]
It follows that the family of $\{v_n\}_{n\in\mathbb{N}}$ is also uniformly bounded. Using the Arzel\`{a}-Ascoli theorem again, $v_n$ converges uniformly to a Lipschitz map  $v$ (up to a subsequence). Meanwhile, it is clear to see that for a given $t\in \R$, $\Phi_{V_n}^t$ converges uniformly to $\Phi_{V}^t$ on each subset of $\mathbb{T}^d\times\R^d$. Note that the rotation vector $\rho$ is non-resonant. According to \cite[Main Result]{FGS}, there exists at most one Lagrangian torus with frequency $\rho$. Then the embedding $w=(u,v)$ is unique. This completes the proof of Proposition \ref{aacoli}.\End

 Theorem \ref{mathe} implies
\[\Delta(\rho,r,\delta):=\overline{W_\rho}\backslash W_\rho\neq \emptyset.\]
Otherwise, we have $W_\rho=\overline{W_\rho}$. It means that $W_\rho$ is both open and closed in $\mathcal{U}^r_\delta$. Then $W_\rho=\mathcal{U}^r_\delta$, which contradicts Theorem \ref{mathe}. By Proposition \ref{aacoli} and Definition \ref{dd1}, for each $V\in \Delta(\rho,r,\delta)$, the flow generated by $H_0+V$ still has a Lagrangian torus with rotation vector $\rho$, which is formulated as a Lipschitz graph $\Psi:\T^d\to \T^d\times\R^d$.  Item (1) and Item (2) of Theorem \ref{excc} are based on the definition of $W_\rho$.

\section{Proof of Theorem \ref{neev11}}\label{77}
The proof is inspired by Herman  \cite[Corollary 3.5]{H1}. Following Herman, we use $\mathrm{Diff}^r_+(\T)$ to denote the group of $C^r$ diffeomorphisms which are isotopy to the identity on $\T$ where $0\leq r\leq \omega$.  Denote $D^r(\T)$ the universal covering space of $\mathrm{Diff}^r_+(\T)$. Given $f\in D^r(\T)$, the rotation number $\rho(f)$ is well-defined.  We consider the following sets:
\[F_\alpha^r:=\{f\in D^r(\T)\ |\ \rho(f)=\alpha\},\]
\[O_\alpha^r:=\{h^{-1}\circ R_\alpha\circ h\ |\ h\in D^r(\T)\},\]
\[H^r:=\{g\in D^r(\T)\ |\ g=Id+\varphi,\ \int_{\T}\varphi(\theta)d\theta=0\}.\]
The Herman-Yoccoz's global theory for circle diffeomorphisms shows $F_\alpha^\infty=O_\alpha^\infty$ if $\alpha$ is a Diophantine number. For simplicity, we still use $f$ to denote an area-preserving twist map defined by
\[f(x,y)=(x+y,y+\varphi(x+y)).\]
An crucial observation by Herman \cite{H1} is that $f$ admits an invariant graph of $\psi$ on $\T$ if and only if there exists $g\in D^r(\T)$ such that
\[\frac{g+g^{-1}}{2}=Id+\frac{1}{2}\varphi.\]
Moreover, $g$ can be given explicitly by $g=Id+\psi$. Define the map $\Phi:F_\alpha^r\to H^r$ via
\[\Phi(g)=\frac{g+g^{-1}}{2}.\]
By \cite[Theorem 2.6.1]{H1}, $\Phi(F_\alpha^\infty)$ is an open subset of $H^\infty$ in the $C^\infty$ topology. 
Let
\[U_\delta:=\{\varphi\in C^\infty(\T)\ |\ \|\varphi\|_{C^{3-\eps}}<\delta\}.\]
It is clear to see $U_\delta$ is open and connected in the $C^\infty$ topology.  Let $W_\delta:=\Phi(F_\alpha^\infty)\cap U_\delta$, and then take the closure of $W_\delta$ denoted by $\overline{W_\delta}$. A similar argument as the proof of Proposition \ref{aacoli} implies $\Phi(F_\alpha^0)$ is closed. It follows that
\[\overline{W_\delta}\subseteq \Phi(F_\alpha^0)\cap U.\]
Moreover, we also know $\overline{W_\delta}\backslash W_\delta$ is not empty. Otherwise, we have $W_\delta=U_\delta$, which contradicts Theorem \ref{mathe} (also \cite[Theorem 4.9]{H1}).

By \cite[Proposition 2.5.4]{H1}, $\Phi$ is injective if $\alpha$ is irrational. Then $\Phi^{-1}$ is well-defined on the image of $\Phi(F_\alpha^\infty)$. Pick $h_*:=Id+\varphi_*\in \overline{W_\delta}\backslash W_\delta$. We have $g_*:=\Phi^{-1}(h_*)$ is not of class $C^\infty$. Otherwise,  $g_*\in F_\alpha^\infty$, which means $h_*\in W_\delta$. Note that the invariant graph $\psi_*$ admitted by $f$ is formulated by $\psi_*=g_*-Id$. Consequently, $\psi_*$ is not of class $C^\infty$.

Take $\delta=\frac{1}{n}$. one can find a sequence $\varphi_n$ such that
\[Id+\varphi_n\in \overline{W_{\frac{1}{n}}}\backslash W_{\frac{1}{n}}.\]
Then we have $\|\varphi_n\|_{C^{3-\eps}}\to 0$ and $\|\varphi_n\|_{C^{3}}\nrightarrow 0$. The fact $\|\varphi_n\|_{C^{3}}\nrightarrow 0$ follows from \cite[Corollary 7.10]{H33}.

 \vspace{2em}

 \noindent\textbf{Acknowledgement.}  Part of this work was done during my visit to Dipartimento di Matematica, Universita degli Studi di Roma `` Tor Vergata" in summer 2023. I would like to thank Prof. A. Sorrentino for his invitation and hospitality. I also warmly appreciate  Prof. A. Bounemoura and Dr. Z. Tong for the discussion on optimality on Item (2) and (3) of Theorem \ref{mathe}. This work is partially under the support of National Natural Science Foundation of China (Grant No. 12122109).

\appendix
\setcounter{table}{0}   
\setcounter{figure}{0}
\setcounter{equation}{0}

\renewcommand{\thetable}{A.\arabic{table}}
\renewcommand{\thefigure}{A.\arabic{figure}}


\addcontentsline{toc}{section}{\sc References}

\begin{thebibliography}{DGOW}
\renewcommand{\baselinestretch}{1.0}
\setlength\itemsep{-1pt}
\small

\bibitem{A} J. Albrecht. {\it On the existence of invariant tori in nearly-integrable Hamiltonian systems with finitely
differentiable perturbations}. Regul. Chaotic Dyn.
\textbf{12} (2007), 281-320.

\bibitem{AK} D. Anosov and A. Katok. {\it New examples in smooth ergodic theory. Ergodic diffeomorphisms}. Trans. Moscow Math. Soc. 23 (1970), 1-35.

\bibitem{Arna0} M. C. Arnaud.  {\it  Fibr\'{e}s de Green et r\'{e}gularit\'{e} des graphes $C^0$-lagrangiens invariants par un
 flot de Tonelli}. Annales Henri Poincare 9 (5) (2008), 881-926.


\bibitem{Arna} M. C. Arnaud.  {\it On a theorem due to Birkhoff}.  Geom. Funct. Anal. {20} (2010), 1307-1316.

\bibitem{Arna2} M. C. Arnaud. {\it A non-differentiable essential irrational invariant curve for a $C^1$ symplectic twist map}. J. Mod. Dyn. 5(3) (2011), 583-591.


\bibitem{ams} M. C. Arnaud, J. Massetti and A. Sorrentino.  {\it On the fragility of periodic tori for families of symplectic twist maps}.  Adv. Math. {429} (2023), Paper No. 109175, 39 pp.

\bibitem{AF} A. Avila and B. Fayad. {\it Non-differentiable irrational curves for  $C^1$  twist map}.
Ergod.
 Th. \& Dynam. Sys. 42 (2022), no. 2, 491-499.

%
%
%
%
%
\bibitem{B} V. Bangert. {\it Mather sets for twist maps and geodesics on
tori}. Dynamics Reported {1} (1988), 1-45.




 \bibitem{BK} D. Bernstein and A. Katok. {\it Birkhoff periodic
orbits for small perturbation of completely integrable Hamiltonian
systems with convex Hamiltonians}. Invent. Math. {88} (1987),
225-241.


\bibitem{B2} U. Bessi.  {\it An analytic counterexample to KAM theorem}.
Ergod. Th. \& Dynam. Sys. {20} (2000), 317-333.




\bibitem{BP} M. Bialy and L. Polterovich. {\it Hamiltonian systems,
Lagrangian tori and Birkhoff's theorem}. Math. Ann. {292}
(1992), 619-627.







%
%

\bibitem{C} C.-Q. Cheng. {\it Non-existence of KAM torus}. Acta Mathmatica
Sinica. {27} (2011), 397-404.




\bibitem{CW}
C.-Q. Cheng and L. Wang. {\it Destruction of Lagrangian torus for
positive definite Hamiltonian systems}. Geom. Funct. Anal. {23} (2013), 848-866.






\bibitem{CI} G. Contreras and R. Iturriaga. {\it  Global minimizers of autonomous Lagrangians}. 22nd  Col\'{o}quio Brasileiro de Matem\'{a}tica. [22nd Brazilian Mathematics Colloquium], Instituto de  Matem\'{a}tica Pura e Aplicada (IMPA), Rio de Janeiro, 1999.



\bibitem{FGS} A. Fathi, A. Giuliani and A. Sorrentino. {\it Uniqueness of invariant Lagrangian graphs in a homology or a cohomology class}.
Ann. Sc. Norm. Super. Pisa Cl. Sci. (5) 8 (2009), no. 4, 659-680.



\bibitem{Fav} J. Favard. {\it Sur les meilleurs proc\'{e}d\'{e}s d'approximation}. B.S.M. 61 (1937), 243-256.




\bibitem{FK} B. Fayad and A. Katok. {\it Constructions in elliptic dynamics}.
Ergod.
 Th. \& Dynam. Sys. 24 (2004), no. 5, 1477-1520.




\bibitem{Fa} B. Fayad. {\it Weak mixing for reparameterized linear flows on the torus}. Ergod.
 Th. \& Dynam. Sys. 22 (2002), no. 1, p. 187-201.



\bibitem{FKO} B. Fayad and R. Krikorian. {\it   Some questions around quasi-periodic dynamics}.
World Scientific Publishing Co. Pte. Ltd., Hackensack, NJ, 2018, 1909-1932.



%
%





 \bibitem{Fo} G. Forni. {\it
Analytic destruction of invariant circles}. Ergod.
 Th. \& Dynam. Sys. {14} (1994), 267-298.

%





%
%
%
%
%
%


 \bibitem{H11} M. R. Herman. {\it Sur la conjugation diff$\acute{e}$rentiable des diff$\acute{e}$omorphismes du cercle $\grave{a}$ des
 rotations}. Publ. Math. IHES {49} (1979), 5-233.

\bibitem{H1} M. R. Herman.  {\it Sur les courbes invariantes par les diff\'{e}omorphismes de
 l'anneau}. Ast\'{e}risque {103-104} (1983), 1-221.
%
\bibitem{H33} M. R. Herman.  {\it Sur les courbes invariantes par les
diff$\acute{e}$omorphismes de
 l'anneau}. Ast$\acute{\text{e}}$risque {144} (1986), 1-243.


\bibitem{H2} M. R. Herman.   {\it In\'{e}galit\'{e}s ``a priori" pour des
tores lagrangiens invariants par des diff\'{e}omrphismes
symplectiques}. Publ. Math. IHES {70}
(1990), 47-101.



%
\bibitem{He6} M. R. Herman.   {\it On the dynamics of
Lagrangian tori invariant by symplectic diffeomorphisms}. In
Progress in variational methods in Hamiltonian systems and elliptic
equations (L'Aquila, 1990). Pitman Res. Notes Math. Ser.
{243} (1992), 92-112.






\bibitem{KH} A. Katok and B. Hasselblatt. {\it Introduction to the Modern Theory of Dynamical Systems} (Encyclopedia of Mathematics and its Applications, 54). Cambridge University Press, Cambridge, 1995.


\bibitem{KO}
Y. Katznelson and D. Ornstein. {\it
Smoothness of invariant curves}.
Proceedings of the Conference in Honor of Jean-Pierre Kahane (Orsay, 1993)
J. Fourier Anal. Appl.,Special Issue(1995), 283-310.








\bibitem{MY} S. Marmi and J.-C. Yoccoz. {\it Some open problems related to small divisors}
Lecture Notes in Math., 1784
Fond. CIME/CIME Found. Subser.
Springer-Verlag, Berlin, 2002, 175-191.







%
 \bibitem{Mm2} J. N. Mather. {\it Non-existence of invariant circles}. Ergod.
 Th. \& Dynam. Sys. {4} (1984), 301-309.
%
%
%

 \bibitem{M4} J. N. Mather. {\it Destruction of invariant circles}. Ergod.
 Th. \& Dynam. Sys. {8} (1988), 199-214.



\bibitem{M5} J. N. Mather. {\it Action minimizing invariant measures for positive definite Lagrangian
systems}. Math. Z. {207} (1991), 169-207.




\bibitem{Mor} M. Morse. {\it A fundamental class of geodesics on any closed surface of genus
greater than one}. Trans. Am. Math. Soc., 26 (1924), 25-60.




%
%
%
%
%
%





\bibitem{NS} V. V. Nemitsky, V. V. Stepanov. {\it Qualitative Theory of Differential Equations}. Princeton, Princeton
University Press, 1960.


\bibitem{NW} P. Ni and L. Wang. {\it Aubry-Mather theory for contact Hamiltonian systems III}. {Sci. China Math.}, to appear.


\bibitem{OP} A. Olvera and N. Petrov. {\it Regularity Properties of Critical Invariant Circles
of Twist Maps, and Their Universality}. SIAM J. Appl. Dyn. Syst.   7 (2008), 962-987.

%

%




\bibitem{P} J. P\"{o}schel.  {\it \"{U}ber invariante Tori in differenzierbaren Hamiltonschen Systemen}. (German)[On invariant tori in differentiable Hamiltonian systems]
Bonner Math. Schriften, 120[Bonn Mathematical Publications]
Universit\"{a}t Bonn, Mathematisches Institut, Bonn, 1980. 103 pp.

\bibitem{P2} J. P\"{o}schel.  {\it KAM below $C^n$}. https://arxiv.org/abs/2104.01866.


%




\bibitem{Sa} D. Salamon. {\it The Kolmogorov-Arnold-Moser theorem}.
Math. Phys. Electron. J. 10 (2004), Paper 3, 37 pp.


\bibitem{SaZ} D. Salamon and E. Zehnder. {\it KAM theory in configuration space}.
Comment. Math. Helv. 64 (1989), no. 1, 84-132.

%

\bibitem{TL} Z. Tong and Y. Li. {\it Towards continuity: Universal frequency-preserving KAM persistence and remaining regularity}.
Communications in Contemporary Mathematics, to appear.

%
%
%
%




\bibitem{W2}L. Wang. {\it Total destruction of Lagrangian tori}. J. Math. Anal. Appl. {410} (2014), 827-836.012).
%


%
\bibitem{W4} L. Wang. {\it Quantitative destruction of invariant circles}. Discrete Contin. Dyn. Syst., 42 (2022), 1569-1583.


%
\bibitem{Z} A. Zygmund. {\it Trigonometric Series}. Third Edition Volumes I \& II combined, with a foreword by Robert Fefferman. Cambridge
University Press, Cambridge, 2002.

\end{thebibliography}
\end{document}